\documentclass[eqno, oneside, draft,11pt]{amsart}

\newcommand{\C}{{\mathbb C}}
\newcommand{\R}{{\mathbb R}}
\newcommand{\Q}{{\mathbb Q}}
\newcommand{\Z}{{\mathbb Z}}
\newcommand{\CP}{{\mathbb CP}}

\usepackage{amsmath}
\usepackage{amsthm}
\usepackage{times}
\usepackage[all]{xy}
\usepackage{amscd}
\usepackage{amsmath,amssymb,latexsym}

\def\Int{\rm Int}

\numberwithin{equation}{section}

\title[Hofer-Zehnder capacity of neighborhoods of symplectic submanifolds]
{Finiteness of the Hofer-Zehnder  capacity of neighborhoods of
symplectic submanifolds}
\author[G. C. Lu]{Guangcun Lu}
\address{Department of Mathematics,\hfill\break\indent  Beijing Normal University, Beijing 100875,   P.R.China}
\email{gclu@bnu.edu.cn} \urladdr{http://math.bnu.edu.cn/\~{}gclu}

\thanks{{\bf IMRN} International Mathematics Research Notices, Volume 2006, Article ID 76520, Pages
1-33.}

\begin{document}

\newtheorem{theorem}{Theorem}[section]
\newtheorem{cor}[theorem]{Corollary}
\newtheorem{corollary}[theorem]{Corollary}
\newtheorem{defn}[theorem]{Definition}
\newtheorem{thm}[theorem]{Theorem}
\newtheorem{definition}[theorem]{Definition}
\newtheorem{example}[theorem]{Example}
\newtheorem{remark}[theorem]{Remark}
\newtheorem{rmk}[theorem]{Remark}
\newtheorem{lemma}[theorem]{Lemma}
\newtheorem{claim}[theorem]{Claim}
\newtheorem{question}[theorem]{Question}
\newtheorem{guess}[theorem]{Conjecture}
\newtheorem{prop}[theorem]{Proposition}
\newtheorem{proposition}[theorem]{Proposition}
\newtheorem{assumption}[theorem]{Assumption}
\newfont{\Bb}{msbm10 scaled\magstephalf}


\begin{abstract}
We use the minimal coupling procedure of Sternberg and Weinstein
and our pseudo-symplectic capacity theory to prove that every
closed symplectic submanifold in any symplectic manifold has an
open neighborhood with finite ($\pi_1$-sensitive) Hofer-Zehnder
symplectic capacity. Consequently,  the Weinstein conjecture holds
near closed symplectic submanifolds in any symplectic manifold.
\end{abstract}

\maketitle

\section{Introduction and results}

The existence of  periodic orbits of Hamiltonian flows near a
closed symplectic submanifold was recently studied by several
authors, (cf. \cite{CiGiKe, GiGu, Ke, Lu2, Schl}). This question
is closely related to the following:

\begin{question}\label{q:1.1}
{\rm   Does every compact symplectic submanifold $B$ in any
symplectic manifold $(M,\omega)$ have a neighborhood with finite
($\pi_1$-sensitive) Hofer-Zehnder symplectic capacity?}
\end{question}

A positive answer to Question~\ref{q:1.1} has some nice
applications in symplectic topology and Hamiltonian dynamics. For
example, it directly implies the existence of (contractible)
Hamiltonian periodic orbits on a generic hypersurface near the
symplectic submanifold. For a closed symplectic submanifold $B$ of
 a geometrically bounded symplectic manifold $(M,\omega)$, if $(M, \omega)$ is
  \emph{symplectically aspherical}, that is,
$\omega|_{\pi_2(M)}=0$ and $c_1(TM)|_{\pi_2(M)}=0$,   Cieliebak,
Ginzburg and  Kerman  \cite{CiGiKe} used symplectic homology to show
that for a smooth function $H:M\to\R$ which attains an isolated
minimum on $B$, the levels $\{H=\epsilon\}$ carry contractible
periodic orbits for a dense set of small values $\epsilon>0$. Under
the same assumptions, it was proved by Ginzburg and G\"urel in
\cite{GiGu} that for a sufficiently small neighborhood $U$ of $B$,
there is a constant $C=C(U)$ such that the Hamiltonian flow of every
smooth function $H$ supported in $U$ with $\min_BH>C$ has a
nontrivial contractible periodic orbit of period less than or equal
to one. In \cite[Corollary 1.3]{Ma}, Macarini showed that for a
closed symplectic submanifold $B$ of any geometrically bounded
symplectic manifold $(M,\omega)$, there exists a neighborhood $U$ of
$B$ in $M$ such that if $H$ is a proper Hamiltonian on $U$ and
constant on $B$, then $H$ has periodic orbits with contractible
projection on $B$ on almost all energy levels. Recently,  Schlenk
\cite{Schl} generalized Hofer's energy-capacity inequality and used
it to prove that for  a closed submanifold $B$ in geometrically
bounded and stably strongly semi-positive symplectic manifolds, if
either $\dim B<{\rm codim}B$, or $\dim B={\rm codim}B$ and $B$ is
not Lagrangian, then $B$ has a small open neighborhood with finite
($\pi_1$-sensitive) Hofer-Zehnder capacity. In \cite{Lu2}, the
author introduced the pseudo symplectic capacity theory and used it
to answer affirmably Question~\ref{q:1.1} for any closed symplectic
submanifold of codimension two in any symplectic manifold.

Very recently, motivated by these and a very useful formula for the
Hofer-Zehnder symplectic capacity of the product of a symplectic
manifold and the standard symplectic ball, (cf. [4, 8, 15, 16, 23,
30]), Kerman \cite{Ke} proposed the following.

\begin{question}\label{q:1.2}
 {\rm Let $B$ be a closed symplectic
submanifold of a symplectic manifold $(M,\omega)$, and let $U_R$
be a symplectic tubular neighborhood of $B$ with (sufficiently
small) radius $R$. Is the Hofer-Zehnder capacity of $U_R$ equal to
$\pi R^2$?}
\end{question}

Here $U_R$ is defined as follows. Let $\pi:(E,\sigma_E)\to B$ be
the symplectic normal bundle of $B$ in $(M,\omega)$, and
$J_E\in{\mathcal J}(E,\sigma_E)$, where ${\mathcal J}(E,\sigma_E)$
is the space of all complex structures on the vector bundle $E\to
B$ compatible with $\sigma_E$, (cf. \cite[page 69]{McSa1}). Denote
by $g_E=\sigma_E\circ(1\times J_E)$ the Hermitian metric on $E$.
Then with a Hermitian connection $\nabla$ on $E$, one can extend a
canonical fiberwise $1$-form $\alpha_E$ on $E\setminus B$ defined
by $\alpha_E(b,v)(\cdot)=g_E(v,v)^{-1}\sigma_E(b)(v,\cdot)$ to a
genuine $1$-form $\alpha$ on $E\setminus B$. In \cite{Ke},
$U_R=\{(b,v)\in E\,|\, g_E(b)(v,v)\le R^2\}$ equipped with the
symplectic form $\omega|_B+ d(g_E(v,v)\alpha)$ was called a {\it
symplectic tubular neighborhood} of $B$.

Clearly, by Weinstein's symplectic neighborhood theorem,
Question~\ref{q:1.2} implies Question~\ref{q:1.1} since $U_R$ is
symplectomorphic to a neighborhood of $B$ in $(M,\omega)$ for
$R>0$ small enough. Using Floer homology, Kerman \cite{Ke}
 affirmably answered Question~\ref{q:1.2} for the closed symplectic submanifold $B$
of dimension $2m$ and codimension $2k$ in a geometrically bounded,
symplectically aspherical manifold $(M,\omega)$ whose unit normal
bundle $S(E)$ is homologically trivial in degree $2m$
 in the sense that $H_{2m}(S(E),\Z_2)=H_{2m}(B,\Z_2)\oplus H_{2(m-k)+1}(B,\Z_2)$.
Actually Questions~\ref{q:1.1}, \ref{q:1.2} may be viewed as
special cases of a general question: \emph{how to compute the
symplectic capacities of a symplectic fibration equipped with a
compatible symplectic form}, which was proposed by C. Viterbo to
the author during his visit to Institut des Hautes \'Etudes
Scientifiques (IHES) in Spring 1999.

The proof of in \cite[Theorem 1.24]{Lu2}, which was partially
motivated by \cite{Bi, Po1},  consists of three steps. The first
step is to construct a suitable projective bundle over the
submanifold and a symplectic form on the total space of the
projective bundle such that the symplectic submanifold has a
neighborhood symplectomorphic to that of the zero section in the
projective bundle. Then we proved that the total space of the
projective bundle is symplectically uniruled with respect to the
chosen symplectic form. Finally, the desired result follows from the
properties of our pseudo symplectic capacity developed in that
paper.  Since we there considered the symplectic submanifold $B$ of
codimension two, its symplectic normal bundle has real rank $2$, and
thus may be viewed a (complex) line bundle. The projectivized bundle
of the sum of the latter and trivial complex line bundle is a
$\CP^1$-bundle, and it is not hard to construct a symplectic form on
the total space of it for which this $\CP^1$-bundle is
symplectically uniruled and naturally contains $B$ as a symplectic
submanifold with the symplectic normal bundle being isomorphic to
that of $B$ in the original symplectic manifold. For the symplectic
submanifolds of higher codimension, the construction of the expected
symplectic form on that kind of projective bundles needs be
elaborated.  This can be completed with the construction of the
minimal coupling by Sternberg and Weinstein. To the knowledge of the
author Polterovich, \cite{Po1} first used the coupling form in the
study of symplectic topology. It had been furthermore used in
\cite{Po2} (also see \cite{Po3}) and \cite{En}.

To state our results, we need to review several notions. Firstly,
for the conveniences of the readers, we recall the definition of
the $\pi_1$-sensitive Hofer-Zehnder capacity $c_{HZ}^\circ$. It
was introduced by the author \cite{Lu1, Lu1+} (denoted by $\bar
C_{HZ}$), and Schwarz \cite{Schw} independently. Recall that in
\cite{HZ1}, a smooth real function $H$ on a symplectic manifold
$(M,\omega)$ is called \emph{admissible}\, if there exist an
nonempty open subset $U$ and a compact subset $K\subset
M\setminus\partial M$ such that
\begin{itemize}
\item[(a)] $H|_U=0$ and $H|_{M\setminus K}=\max H$;

\item[(b)] $0\le H\le\max H$;

\item[(c)] $\dot x=X_H(x)$ has no nonconstant fast periodic
solutions.
\end{itemize}
 Here $X_H$ is defined by
$i_{X_H}\omega=dH$, and ``fast'' means ``of period less than
$1$''. Denote by  ${\mathcal H}_{ad}(M,\omega)$ the set of
admissible Hamiltonians on $(M,\omega)$. The Hofer-Zehnder
symplectic capacity  of $(M,\omega)$ is defined by
$c_{HZ}(M,\omega) = \sup \left\{\max H\,|\, H\in {\mathcal
H}_{ad}(M,\omega) \right\}$. If the condition (c) in the
definition of the admissibility is replaced by the condition
$(c)^\circ$,\;$\dot x=X_H(x)$ has no nonconstant fast and
contractible  periodic solutions, the corresponding function $H$
is said to be $0$-\emph{admissible}. Let ${\mathcal
H}^\circ_{ad}(M,\omega)$ be the set of $0$-admissible Hamiltonians
on $(M,\omega)$. Then the $\pi_1$-sensitive Hofer-Zehnder capacity
of $(M,\omega)$ is defined by
\begin{equation}\label{e:1.0}
c^\circ_{HZ}(M,\omega) = \sup \left\{\max H\,|\, H\in {\mathcal
H}^\circ_{ad}(M,\omega) \right\}.
\end{equation}
It always holds that ${\mathcal W}_G\le c_{HZ}(M,\omega)\le
c_{HZ}^\circ(M,\omega)$ for them and the Gromov symplectic width
 ${\mathcal W}_G$.

Recall that a \emph{symplectic fibration} $\Pi:M\to B$ with
symplectic fibre $(F,\sigma)$ is a fibration whose structure group
is a subgroup of ${\rm Symp}(F,\sigma)$. In this case, for any
local trivialization $\Phi:\Pi^{-1}(U)\to U\times F$ and any $b\in
U$, there exists a natural symplectic form
$\sigma_b=\Phi_b^\ast\sigma$ on the fiber $F_b$ which is
independent of the choice of the local trivialization $\Phi$,
where $\Phi_b:F_b\to F$ is the restriction of $\Phi$ to $F_b$
followed by the projection onto $F$. A symplectic form $\omega$ on
$M$ is said to be \emph{compatible} with the symplectic fibration
$\Pi$ if $\sigma_b=\iota_b^\ast\omega$ for each inclusion
$\iota_b:F_b\hookrightarrow M$ of the fibre, that is, each fibre
$(F_b,\sigma_b)$ is a symplectic submanifold of $(M,\omega)$.

Let $\pi:(E,\bar\omega)\to B$ be a $2n$-dimensional symplectic
vector bundle over a compact symplectic manifold $(B, \beta)$
(with or without boundary). It may be naturally viewed as a
symplectic fibration with symplectic fibre
$(\R^{2n},\omega_{std})$, and hence each fibre $E_b$ carries a
natural symplectic structure $(\omega_{\rm std})_b$, (cf. Section
3 for details). Throughout this paper, we use $\omega_{std}$ (or
$\omega^{(n)}_{std}$ if necessary) to denote the standard
symplectic form on $\R^{2n}$. With a $\bar\omega$-compatible
complex structure $J_E\in{\mathcal J}(E, \bar\omega)$, one gets a
Hermitian structure $(\bar\omega, J_E, g_{J_E})$ on $E$. Denote by
${\mathcal F}:{\rm U}(E)\to B$  the bundle of unitary frames of
$E$. It is a principal $U(n)$-bundle, and $E$ is an associated
bundle, $E={\rm U}(E)\times_{U(n)}\C^{n}$, (see (\ref{e:3.11}) for
an explicit identification). Hereafter, the unitary group $U(n)$
acts on $\C^n$ via
\begin{equation}\label{e:1.1}
U\cdot (z_1,\cdots, z_n)=(z_1,\cdots,z_n)U\quad\forall U=X+ iY\in
U(n).
\end{equation}
This action is Hamiltonian with respect to $\omega_{\rm std}$, and
 has the moment map
\begin{equation}\label{e:1.2}
\mu_{U(n)}:\C^n\to {\frak u}(n),\;z=(z_1,\cdots,z_n)\mapsto
\frac{i}{2}z^t\bar z
\end{equation}
after identifying the Lie algebra ${\frak u}(n)=T_{I_n}U(n)$ with
its dual ${\frak u}(n)^\ast$  via the inner product $(X,Y)={\rm
Tr}(\bar X^tY)$.  Denote by
\begin{equation}\label{e:1.3}
{\rm D}_\varepsilon(E)={\rm U}(E)\times_{U(n)} B^{2n}(\varepsilon)
\end{equation}
the open disk bundle of radius $\varepsilon>0$.  Here
$B^{2n}(\epsilon)=\{{\bf z}\in\C^n\,|\, |{\bf z}|<\varepsilon\}$.
Let ${\mathcal A}({\rm U}(E))$ be the affine space of all connection
(${\frak u}(n)$-value) $1$-forms on ${\rm U}(E)$. By the
construction due to Sternberg and Weinstein, (cf.
Theorem~\ref{th:3.2}($8^\circ$)), for each $A\in {\mathcal A}({\rm
U}(E))$ there are  $0<\varepsilon_0=\varepsilon_0({\mathcal F},
A)\le 1$ and a canonical symplectic form $\bar\omega_A$ in ${\rm
D}_{\varepsilon_0}(E)$ such that
\begin{equation}\label{e:1.4}
\bar\omega_A|_{{\rm D}_{\varepsilon_0}(E)_b}=(\omega_{\rm
std})_b|_{{\rm D}_{\varepsilon_0}(E)_b}\quad\forall b\in B,
\end{equation}
and that the symplectic normal bundle of the zero section $0_E$ in
$({\rm D}_{\varepsilon_0}(E), \bar\omega_A)$ is $(E,\bar\omega)$.
 For each $t>0$ let
\begin{equation}\label{e:1.5}
\bar\omega_A^t:=\pi^\ast\beta+ t(\bar\omega_A-\pi^\ast\beta).
\end{equation}
It restricts to $t(\omega_{\rm std})_b$ on each fibre ${\rm
D}_{\varepsilon_0}(E)_b$, and for every
$0<\varepsilon<\varepsilon_0$ there is $t_0=t_0(A,\varepsilon)>0$
such that $\bar\omega_A^t$ is also symplectic in ${\rm
D}_\varepsilon(E)$ for each $0<t<t_0(A,\varepsilon)$, the reader
may refer to Theorems 3.3(v). The pair $({\rm D}_{\varepsilon}(E),
\bar\omega^t_A)$ is our \emph{symplectic tubular neighborhood} of
$B$. Our main result is the following theorem.

\begin{theorem}\label{th:1.3}
 Let  $B$ be any closed symplectic submanifold in any symplectic
manifold $(M, \omega)$ and let $\pi:(E,\bar\omega)\to B$ be the
symplectic normal bundle of $B$ in $(M,\omega)$. For a given
$J_E\in{\mathcal J}(E,\bar\omega)$,  let ${\mathcal F}:{\rm U}(E)\to
B$ be the corresponding bundle of unitary frames. Then for a
connection $1$-form $A\in {\mathcal A}({\rm U}(E))$ and
$\bar\omega_A^t$ in (\ref{e:1.5}) with $\beta=\omega|_B$,  it holds
that
\begin{equation}\label{e:1.6}
c_{HZ}({\rm D}_{\varepsilon}(E),\bar\omega_A^t)\le
c_{HZ}^\circ({\rm D}_{\varepsilon}(E),\bar\omega_A^t)\le\pi
t\varepsilon^2
\end{equation}
for any $0<\varepsilon<\varepsilon_0({\mathcal F}, A)$ and
$0<t<t_0(\varepsilon, A)$.  Furthermore
\begin{equation}\label{e:1.7}
{\mathcal W}_G({\rm
D}_{\varepsilon}(E),\bar\omega_A^t)=c_{HZ}({\rm
D}_{\varepsilon}(E),\bar\omega_A^t)=c_{HZ}^\circ({\rm
D}_{\varepsilon}(E),\bar\omega_A^t)=\pi t\varepsilon^2
\end{equation}
for sufficiently small $\varepsilon>0$ and $t>0$. So for any
$\epsilon>0$ there exists an open neighborhood $W$ of $B$ in $M$
such that
\begin{equation}\label{e:1.8}
{\mathcal
W}_G(W,\omega)=c_{HZ}(W,\omega)=c^\circ_{HZ}(W,\omega)<\epsilon.
\end{equation}
\end{theorem}

From the proof of Theorem~\ref{th:1.3}, it is easily seen that for
the above \linebreak $(U_R, \omega|_B+ d(g_E(v,v)\alpha))$ there
exist $0<\epsilon\ll\varepsilon$ and a symplectic embedding from
$({\rm D}_{\epsilon}(E), \bar\omega^t_A)$ into $(U_R, \omega|_B+
d(g_E(v,v)\alpha))$ which maps the zero section $0_E$ onto $0_E$.
Conversely, there exist $0<r\ll R$ and a symplectic embedding from
$(U_r, \omega|_B+ d(g_E(v,v)\alpha))$ into $({\rm
D}_{\varepsilon}(E), \bar\omega^t_A)$ mapping the zero section
$0_E$ onto $0_E$.

It is well known that the Hofer-Zehnder capacity and pseudo
symplectic capacity are closely related to the famous Weinstein
conjecture in \cite{We}, (cf. \cite{HZ1, Lu2}). Every compact smooth
hypersurface ${\mathcal S}$ in a symplectic manifold $(M, \omega)$
determines a distinguished line bundle $T{\mathcal S}\supset
{\mathcal L}_{\mathcal S}\rightarrow {\mathcal S}$ whose fiber at
$x\in {\mathcal S}$ is given by $T_x{\mathcal S}\cap(T_x{\mathcal
S})^\omega$. A closed characteristic  of ${\mathcal S}$ is an
embedded circle $P\subset {\mathcal S}$ satisfying $TP={\mathcal
L}_{\mathcal S}|_P$. The hypersurface ${\mathcal S}$ is said to be
of \emph{contact type} if there exists a Liouville vector field $X$
(i.e., $L_X\omega=\omega$) in some neighborhood of it which is
transversal to ${\mathcal S}$ everywhere. Weinstein \cite{We}
conjectured  that \emph{every hypersurface ${\mathcal S}$ of contact
type in symplectic manifolds carries a closed characteristic}. After
Viterbo \cite{V} first proved it in the standard Euclidean
symplectic space this conjecture was proved in many symplectic
manifolds, (cf. \cite{Gi, Lu2, Schl} and the references therein for
a detailed description of the progress on this question). As a
direct consequence of (\ref{e:1.8}), we get the following corollary.

\begin{corollary}  \label{cor:1.4}
Let $W$ be a neighborhood of $B$ as in (\ref{e:1.8}). Then for every
smooth function $H:M\to\R$ supported in $W$ and with $\max H-\min
H>c_{HZ}^\circ(W)$, its Hamiltonian flow  has a nontrivial
contractible periodic orbit of period less than or equal to one.
Moreover, for a compact hypersurface ${\mathcal S}$ contained in $W$
and every thickening of ${\mathcal S}$ in $W$, $\psi:{\mathcal
S}\times (-1, 1)\to W$ there is a closed characteristic on
${\mathcal S}_t:=\psi({\mathcal S}\times t)$ for almost every $t\in
(-1, 1)$. In particular, the Weinstein conjecture holds near a
closed symplectic submanifold $B$ in any symplectic manifold $(M,
\omega)$.
\end{corollary}

So Corollary~\ref{cor:1.4} generalizes the corresponding results
in \cite{CiGiKe, GiGu, Ke, Lu2, Ma}. The proof from (\ref{e:1.8})
to Corollary~\ref{cor:1.4} is standard, see \cite{HZ2, MaSchl,
St}. Another direct consequence of (\ref{e:1.8}) is the following
corollary.

\begin{corollary}  \label{cor:1.5}
Let $(N,\sigma)$ be any closed symplectic manifold. Then there is
an open neighborhood $U$ of the zero section in the twisted
cotangent bundle $(T^\ast N, \omega_{\rm can}+ \pi^\ast\sigma)$
such that $c_{HZ}^\circ(U, \omega_{\rm can}+
\pi^\ast\sigma)<\infty$.
\end{corollary}

This implies the existence of contractible periodic orbits of a
charge on the symplectic manifold $N$ subject to the magnetic
field $\sigma$ on almost every sufficiently small energy level. A
more general result was recently obtained in \cite{Schl} by a
different method. The reader may also find some related
interesting results in \cite{Ke}. Equation (\ref{e:1.8}) also
implies some nonsqueezing phenomenon in symplectic geometry which
was first discovered by Gromov in his celebrated paper \cite{Gr}.
For example, it implies that a given standard symplectic ball
$(B^{2n}(r),\omega_{\rm std})$ can not be symplectically embedded
in a small neighborhood of a closed symplectic submanifold in any
symplectic manifold of dimension $2n$.

The organization of the paper is as follows. In the next section we
first review the minimal coupling procedure by Sternberg and
Weinstein following [37, Appendix], and then prove the main result
Theorem~\ref{th:2.5} therein. In Section 3, we shall construct the
symplectic forms on the total space of the projectivized bundle of
the sum of the symplectic normal bundle and the trivial line bundle,
and study their properties. Theorems~\ref{th:3.2}, \ref{th:3.3}
summarize our main results in that section. The main result
Theorem~\ref{th:1.3} will be proved in Section 4 after main results
Theorems~\ref{th:4.4}-\ref{th:4.5} in that section. \vspace{2mm}

\section{Minimal coupling and symplectic reduction}
\setcounter{equation}{0}

In this section, we first briefly review how to use the
 minimal coupling form procedure of Sternberg
and Weinstein to construct the symplectic structures on the
associated bundle,  and then point out some properties of such
symplectic structures. Those properties are needed in our
arguments and  are easily proved. Our main reference is
\cite[Appendix]{SjLe}.

 Let $G$ be a compact Lie group with Lie algebra ${\frak
g}$, and let $\langle\cdot,\cdot\rangle$ denote the pairing between
the dual Lie algebra ${\frak  g}^\ast$ and ${\frak g}$. For a
Hamiltonian $G$-space $(M,\omega)$, let $\mu_G:M\to {\frak g}^\ast$
denote the ${\rm Ad}^\ast$-equivariant moment map. The famous
Marsden-Weinstein reduction theorem is shown as follows.

\begin{theorem}\label{th:2.1}
Suppose that $0$ is a regular value of $\mu_G$ and that the group
$G$ acts freely and properly on $\mu_G^{-1}(0)$. Then
$\mu_G^{-1}(0)$ is a coisotropic submanifold in $(M,\omega)$ and
the corresponding isotropic foliation is given by the orbits of
$G$. Moreover, there exists a unique symplectic form $\omega_G$ on
the Marsden-Weinstein quotient $M/\!\!/G:=\mu_G^{-1}(0)/G$ such
that $q_G^\ast\omega_G=\omega|_{\mu_G^{-1}(0)}$, where
$q_G:\mu_G^{-1}(0)\to M/\!\!/G$ is the quotient projection.
\end{theorem}

 A symplectic submanifold
$N\subset M$ is called a \emph{Hamiltonian $G$-subspace} of
$(M,\omega)$ if $g\cdot N=\{g\cdot x\,|\, x\in N\}\subset N$ for
each $g\in G$. Clearly, for a Hamiltonian $G$-subspace  $N\subset
M$ of $(M,\omega)$, the moment map of the $G$-action on
$(M,\omega)$ restricts to that of the induced $G$-action on
$(N,\omega|_N)$. By Theorem~\ref{th:2.1}, one easily gets the
following proposition.

\begin{proposition}\label{prop:2.2}
For a Hamiltonian $G$-subspace $N$ of $(M,\omega)$, if $0$ is also
a regular value of $\mu_G|_N$ and $G$ also acts freely and
properly on $(\mu_G|_N)^{-1}(0)$, then the restriction of
$\omega_G$ to $N/\!\!/G=(\mu_G|_N)^{-1}(0)/G$ is exactly the
unique symplectic structure from symplectic reduction of the
Hamiltonian $G$-subspace $(N, \omega|_N)$.
\end{proposition}

Now let $\pi_P:P\to B$ be a principal $G$-bundle over a compact
symplectic manifold $(B,\beta)$, and let $\pi_2:P\times{\frak
g}^\ast\to {\frak  g}^\ast$ be the natural projection to the second
factor. Every connection $1$-form $A\in{\mathcal A}(P)$ yields a
corresponding \emph{minimal coupling form} on $P\times{\frak
g}^\ast$ defined by
\begin{equation}\label{e:2.1}
\delta_A={\pi_P}^\ast\beta-d\langle \pi_2, A\rangle.
\end{equation}
By the transformation properties of connections, $\langle \pi_2,
A\rangle$, and hence $\delta_A$, is $G$-invariant under the diagonal
action
\begin{equation}\label{e:2.2}
 g\cdot (p, \zeta)=(g^{-1}\cdot p, {\rm Ad}(g)^\ast\zeta).
\end{equation}
 Moreover, the injection
$\iota:P\to P\times{\frak  g}^\ast$ given by $\iota(p)=(p,0)$ is
$G$-equivariant and $\iota^\ast\delta_A=\pi_P^\ast\beta$. It is
also easily checked that $\delta_A$ is nondegenerate at all points
of $\iota(P)$.

\begin{theorem}\label{th:2.3}
Under the above assumptions, there exists a neighborhood
${\mathcal W}$ of zero in ${\frak  g}^\ast$ (which depends on $A$
and on the principal bundle $\pi_P$, and  can also be required to
be invariant under the coadjoint action of $G$) such that the
following hold.
\begin{itemize}
\item[(i)] The  minimal coupling form $\delta_A$ restricts to a
symplectic form on $P\times{\mathcal W}$.

 \item[(ii)] The action of $G$ on $P\times{\frak  g}^\ast$
 given by (\ref{e:2.2})
  is Hamiltonian with respect to $\delta_A$ on $P\times{\mathcal W}$
 and has the moment map, $-\pi_2:P\times{\frak  g}^\ast\to {\frak
 g}^\ast$.

 \item[(iii)] For two different connection $1$-forms $A_1$ and
 $A_2$ on $P$, let ${\mathcal W}_1$ (resp., ${\mathcal W}_2$) be the corresponding neighborhood
of zero in ${\frak  g}^\ast$ such that the minimal
 coupling form $\delta_{A_1}$ (resp., $\delta_{A_2}$) is symplectic
   on $P\times{\mathcal W}_1$ (resp., on $P\times{\mathcal W}_2$). Then there are {\rm
   smaller} neighborhoods
of zero in ${\frak  g}^\ast$, ${\mathcal
W}^\star_1\subset{\mathcal W}_1$ and ${\mathcal
W}_1^\star\subset{\mathcal W}_2$ such that there exists a
symplectomorphism
 from $(P\times{\mathcal W}^\star_1, \delta_{A_1})$ onto $(P\times{\mathcal W}^\star_2,
 \delta_{A_2})$ that not only commutes with the action of $G$ but also
 restricts to the identity on $P\times\{0\}$.
  \end{itemize}
\end{theorem}

\begin{theorem}\label{th:2.4}
Under the assumptions of Theorem~\ref{th:2.3}, furthermore assume
that $(F,\sigma)$ is a Hamiltonian $G$-space with moment map
$\mu_G^F:F\to {\frak  g}^\ast$ satisfying
\begin{equation}\label{e:2.3}
\mu_G^F(F)\subset{\mathcal W}.
\end{equation}
Then the diagonal action of $G$ on $W:=(P\times{\mathcal W})\times
F$ is Hamiltonian and the corresponding moment map is given by
\begin{equation}\label{e:2.3.1}
\mu_G^W: W\to {\frak  g}^\ast,\;(p, g^\ast, f)\mapsto
\mu_G^F(f)-g^\ast.
\end{equation}
Moreover, there exists a symplectic form $\omega_A$ on the total
space of the associated fibre bundle $\pi_F:M:=P\times_GF\to B$
such that the map $(F,\sigma)\hookrightarrow (M, \omega_A)$ is a
symplectic embedding, that is, $\omega_A|_{M_b}=\sigma_b$ for any
$b\in B$, where $\sigma_b$ is the  symplectic form on $M_b$
obtained from the symplectic fibration $M\to B$ with symplectic
fibre $(F,\sigma)$.  Consequently, each fibre is a symplectic
submanifold.
\end{theorem}

Indeed, it is easy to check that $0$ is a regular value of
$\mu_G^W$, and that
\begin{equation}\label{e:2.3.2}
\widehat\Upsilon: P\times F\to (\mu_G^W)^{-1}(0),\;(p, f)\mapsto
(p,\mu_G^F(f), f)
\end{equation}
is a $G$-equivariant diffeomorphism. Moreover, $G$ acts freely on
$(\mu_G^W)^{-1}(0)$. The Marsden-Weinstein reduction procedure
yields a unique symplectic form $\omega'_A$ on the quotient
\linebreak $(\mu_G^W)^{-1}(0)/G$ whose pullback under the quotient
projection \linebreak $\Pi_W:(\mu_G^W)^{-1}(0)\to
(\mu_G^W)^{-1}(0)/G $ is equal to the restriction of
$\delta_A\oplus\sigma$ to $(\mu_G^W)^{-1}(0)$. Let $\Upsilon:
P\times_GF\to (\mu_G^W)^{-1}(0)/G$ be the diffeomorphism induced by
$\widehat\Upsilon$. Then the symplectic form
$\omega_A=\Upsilon^\ast\omega'_A$ satisfies the desired
requirements. \emph{Note that} shrinking ${\mathcal W}$ while
preserving (\ref{e:2.3}), one obtains the same symplectic form
$\omega_A$ on $P\times_GF$. The following theorem summarizes some
related properties of the above constructions. Their proofs are
easy, (cf. \cite{Lu7}).

\begin{theorem}\label{th:2.5}
Under the assumptions of Theorem~\ref{th:2.4}, the following
properties hold.
\begin{itemize}
 \item[(i)] Let $(F^\sharp,\sigma^\sharp)$ be another
Hamiltonian $G$-space of dimension $\dim F$ with moment map
satisfying (\ref{e:2.3}). If there exists a symplectic embedding
$\varphi:(F^\sharp,\sigma^\sharp)\to (F,\sigma)$ which commutes
with the Hamiltonian actions of $G$ on $(F^\sharp, \sigma^\sharp)$
and $(F,\sigma)$, that is, $g\cdot\varphi(x)=\varphi(g\cdot x)$
for all $g\in G$ and $x\in F^\sharp$, then  the bundle embedding
$\varphi_P:M^\sharp:=P\times_GF^\sharp\to M$ induced by $\varphi$
is a symplectic embedding from $(M^\sharp,\omega^\sharp_A)$ to
$(M,\omega_A)$, where $\omega^\sharp_A$ is the symplectic form on
$M^\sharp$ constructed as above.
 Furthermore, if $F_0$ is a Hamiltonian $G$-subspace of
$(F,\sigma)$, then $P\times_GF_0$ is a symplectic submanifold in
$(M,\omega_A)$.

\item[(ii)] If $F$ is a vector space and $\mu_G^F(0)=0$, then the
zero section $Z_0:=P\times_G\{0\}\subset M$ is a symplectic
submanifold. More precisely,
$\omega_A|_{Z_0}=\pi_F^\ast\beta|_{Z_0}$. Consequently, $Z_0$ is a
symplectic submanifold in $(M,\omega_A)$ and the symplectic normal
space $(T_{(b,0)}Z_0)^{\omega_A}$ at any point $(b,0)\in Z_0$ is
exactly the symplectic vector space $(M_b,\sigma_b(0))$ (because
$T_{(b,0)}M=T_{(b,0)}Z_0\oplus T_0M_b=T_{(b,0)}Z_0\oplus M_b$).
 \item[(iii)]  For any compact symplectic submanifold
$B^\circ\subset B$, (since the connection form $A$ can always
restrict to a connection form on the restriction principal bundle
$P^\circ:=P|_{B^\circ}$ ),  the corresponding minimal coupling form
$\delta_A^\circ=\pi_{P^\circ}^\ast\beta-d\langle\pi_2,A|_{P^\circ}\rangle$
on $P^\circ\times{\frak  g}^\ast$ is equal to the restriction of the
coupling form $\delta_A$ on $P\times{\frak  g}^\ast$ to
$P^\circ\times{\frak  g}^\ast$. By shrinking ${\mathcal W}$, assume
that $\delta_A$ (resp., $\delta_A^\circ$) is nondegenerate on
$P\times{\mathcal W}$ (resp., $P^\circ\times{\mathcal W}$). Then
 $M^\circ:=(P|_{B^\circ})\times_GF$ is also a
symplectic submanifold in $(M,\omega_A)$ and $\omega_A|_{M^\circ}$
is exactly equal to the symplectic form constructed in the above
method from the restriction of the connection form $A$ on the
principal $G$ subbundle $P^\circ\to B^\circ$.

  \item[(iv)] In the trivial principal bundle $P=B\times G$, there
  exists a canonical flat connection defined by ${\mathcal
  H}_u^{cf}={\rm Ker}(\Pi_{2\ast}):T_u(B\times G)\to T_gG$, where
  $\Pi_2:B\times G\to G$ is the projection on the second factor.
The corresponding connection form is given by $A_{\rm
can}:=\Pi_2^\ast\theta$. Here $\theta$ is the canonical left
invariant  ${\frak  g}$-valued $1$-form defined by
$\theta(a)(\widetilde X(a))=X\in{\frak  g}$ for $a\in G$,
$X\in{\frak g}$, where $\widetilde X$ is the unique left invariant
vector field on $G$ which has value $X$ at $e$.
 Then the symplectic form $\omega_{A_{\rm can}}$ on
 $P\times_GF=B\times F$ is equal to $\beta\oplus\sigma$.
 Here the condition (\ref{e:2.3}), of course, has been assumed,
 but the present ${\mathcal W}$ depends merely on $G$ itself.
 \item[(v)] Let $A_1$ and $A_2$ be two connection forms on $P$.
 Suppose that $\mu_F^G(F)$ is contained in the intersection of the
 open subsets ${\mathcal W}^\star_1$ and ${\mathcal W}^\star_2$ in
 Theorem~\ref{th:2.3}(iii). Then there exists a bundle isomorphism
$\Phi: P\times_GF\to P\times_GF$ which sends $\omega_{A_1}$ to
$\omega_{A_2}$, i.e., $\Psi^\ast\omega_{A_2}=\omega_{A_1}$.
\end{itemize}
\end{theorem}

\begin{remark}\label{rm:2.6}
{\rm  For a  Hamiltonian $G$-space  $(F,\sigma)$ with moment map
$\mu_G^F:F\to {\frak  g}^\ast$, there always exists
$\varepsilon_0=\varepsilon_0({\mathcal W}, F,\sigma, G)$ such that
$\varepsilon^2\mu_G^F(F)\subset{\mathcal W}$ for any
$\varepsilon\in (0, \varepsilon_0]$. Since the $G$-action  is also
Hamiltonian with respect to $\varepsilon^2\sigma$ and  the
corresponding moment $\mu_G^\varepsilon=\varepsilon^2\mu_G:F\to
{\frak g}^\ast$, one can always obtain a family of deformedly
equivalent symplectic forms
$\{\omega_{A\varepsilon}\,|\,0<\varepsilon\le\varepsilon_0\}$ on
$P\times_GF$. Here the reason that we use $\varepsilon^2\sigma$
(and thus $\varepsilon^2\mu_G^F$), instead of $\varepsilon\sigma$
(and $\varepsilon\mu_G^F$), will be seen in next section. }
\end{remark}

\section{Symplectic forms on projective bundles}
\setcounter{equation}{0}

In this section, we shall construct two families of symplectic
forms on the disk bundle and projective bundle, and give their
properties. Firstly,  a $2n$-dimensional symplectic vector bundle
$\pi:(E,\bar\omega)\to B$ \emph{may be naturally viewed as a
symplectic fibration with symplectic fibre
$(\R^{2n},\omega_{std})$}.  Let $\bar\omega_{\rm std}$
 denote the standard skew-symmetric bilinear
map on $\R^{2n}$.
 Taking any $J_E\in{\mathcal
J}(E,\bar\omega)$ and setting  $g_{J_E}:E\times
E\to\R,\;g_{J_E}(u,v)=\bar\omega(u,J_Ev)$, one gets a Hermitian
structure $(\bar\omega, J_E,g_{J_E})$ on $E$, (cf. \cite{McSa1}).
Then one can choose an open cover $\{U_\alpha\}_{\alpha\in\Lambda}$
of $B$ such that for each $\alpha\in\Lambda$, there exists a unitary
trivialization
\begin{equation}\label{e:3.0}
U_{\alpha}\times\R^{2n}\to
E|_{U_\alpha}:(b,v)\mapsto \Phi_{\alpha}(b)v
\end{equation}
satisfying $\Phi_{\alpha}^\ast J_E=J_{\rm std}$,
$\Phi_{\alpha}^\ast\bar\omega=\bar\omega_{\rm std}$, and
$\Phi_\alpha^\ast g_J=g_{\rm std}$.
 Here $g_{\rm std}$ and $J_{\rm std}$ are the standard inner
product and complex structure on $\R^{2n}$,  respectively. As a
symplectic fibration $E\to B$, each fibre $E_b$ carries a natural
symplectic structure
\begin{equation}\label{e:3.1}
(\omega_{\rm std})_b:=(\Phi_\alpha(b)^{-1})^\ast\omega_{\rm
std}\;\;({\rm if}\; b\in U_\alpha)
\end{equation}
to satisfy $(\omega_{\rm std})_b|_{T_{0}E_b}=\bar\omega_b$.

The Lie algebra ${\frak  u}(n)=\{X\in M_n(\C)\,|\, \bar X^t=-X\}$
of $U(n)$ carries an invariant (real) inner product defined by
$(X,Y)={\rm Tr}(\bar X^tY)$.   By Riesz theorem, for each
$f\in{\frak  u}(n)^\ast$, there exists a unique $\zeta_f\in{\frak
u}(n)$ such that $\langle f,\xi\rangle=(\xi,\zeta_f)$ for any
$\xi\in{\frak  u}(n)$. Since $ \zeta_{Ad_g^\ast
f}=Ad_{g^{-1}}\zeta_f$ for any $f\in{\frak u}(n)$, and
\begin{equation}\label{e:3.2}
{\frak  u}(n)^\ast\to{\frak  u}(n), f\mapsto\zeta_f
\end{equation}
is a real vector space isomorphism,  we may identify ${\frak
u}(n)^\ast$ with ${\frak  u}(n)$ and $Ad_g^\ast$ with
$Ad_{g^{-1}}$. It is under such identifications that the moment
map of the action of $U(n)$ on $(\C^n,\omega_{\rm std})$ in
(\ref{e:1.1}) is given by (\ref{e:1.2}).

Let $J_{\rm FS}$ be the standard complex structure
 on the $n$-dimensional complex projective space
$\CP^n=(\C^{n+1}-\{0\})/\C^\ast$.
 The standard action of ${\rm U}(n+1)$ on $\C^{n+1}$ as in (\ref{e:1.1})
induces a natural one on $\CP^n$:
\begin{equation}\label{e:3.3}
U\cdot [z_0,\cdots, z_n]=[(z_0,\cdots,z_n)U]
\end{equation}
for any $[z_0,\cdots,z_n]\in\CP^n$ and $U\in U(n+1)$. Let
$\omega_{\rm FS}$ be the unique ${\rm U}(n+1)$-invariant K\"ahler
form on $\CP^n$ associated with the \emph{Fubini-Study metric}
which has integration $\pi$ on $\CP^1\subset\CP^n$. The action in
(\ref{e:3.3}) is Hamiltonian with respect to $\omega_{\rm FS}$ ,
and in homogeneous coordinates $[z_0,\cdots, z_n]$, the
corresponding moment map $\mu'_{U(n+1)}:\CP^n\to {\frak
u}(n+1)^\ast$ is given by
\begin{equation}\label{e:3.3.1}
\langle\mu'_{U(n+1)}([z]),\zeta\rangle=\frac{i}{2}\frac{\sum_{j,k}\zeta_{jk}\bar
z_j z_k}{\sum_j|z_j|^2}.
\end{equation}
Identifying ${\frak  u}(n+1)^\ast$ with ${\frak  u}(n+1)$ via the
isomorphism as in (\ref{e:3.2}),  the moment map
$\mu'_{U(n+1)}:\CP^{n+1}\to {\frak  u}(n+1)$ is written as
\begin{equation}\label{e:3.4}
\mu'_{U(n+1)}([z])=\frac{i}{2}\frac{z^t\bar
z}{\sum_j|z_j|^2}=\biggr(-\frac{1}{2}\frac{(\bar z_j z_k-\bar z_k
z_j)}{\sum_j|z_j|^2}\biggl).
\end{equation}
Consider the Lie group inclusion homomorphism
\begin{equation}\label{e:3.4.1}
\varphi:U(n)\to U(n+1),\;U\mapsto \left(\begin{array}{cr}
U& 0\\
0& 1\end{array}\right),
\end{equation}
we easily get the following proposition.

\begin{proposition}\label{prop:3.1}
The action of $U(n)$ on $\CP^n$ via
\begin{equation}\label{e:3.5}
U\cdot [z_0,\cdots,z_n]=[(z_0,\cdots,z_{n-1})U, z_n]
\end{equation}
is Hamiltonian and the moment map $\mu''_{U(n)}:\CP^n\to {\frak
u}(n)$ is given by
\begin{equation}\label{e:3.5.1}
\mu''_{U(n)}:\CP^{n}\to {\frak  u}(n),\;z\mapsto
\frac{i}{2}\frac{Z^t\bar Z}{\sum^n_{l=0}|z_l|^2}.
\end{equation}
 Here $Z=(z_0,\cdots, z_{n-1})$ and  ${\frak u}(n)^\ast$ has been identified with ${\frak u}(n)$
via (\ref{e:3.2}).
\end{proposition}

 Using the action in (\ref{e:3.5}), one can form a $\CP^n$-bundle
\begin{equation}\label{e:3.6}
{\mathcal P}:{\rm U}(E)\times_{U(n)}\CP^{n}\to B,
\end{equation}
which is exactly the projective bundle ${\rm P}(E\oplus\C)$.
Denote by
\begin{equation}\label{e:3.7}
Z_0:={\rm P}(\{0\}\oplus \C)\quad{\rm and}\quad Z_\infty:={\rm
P}(E\oplus\{0\})
\end{equation}
the zero section and divisor at infinity respectively. The
symplectic  embedding
\begin{equation}\label{e:3.8}
\varphi:(B^{2n}(1),\omega_{\rm std})\to (\CP^n,\omega_{\rm
FS}),\;{\bf z}\mapsto [{\bf z}, \sqrt{1-|{\bf z}|^2}]
\end{equation}
satisfies
\begin{equation}\label{e:3.9}
\varphi(U\cdot{\bf z})=[U\cdot{\bf z},\sqrt{1-|{\bf
z}|^2}]=U\cdot\varphi({\bf z}),
\end{equation}
and hence induces an embedding
\begin{equation}\label{e:3.10}
\varphi_{\mathcal P}:{\rm U}(E)\times_{U(n)}B^{2n}(1)\to {\rm
U}(E)\times_{U(n)}\CP^{n},\; [f, {\bf z}]\mapsto [f,\varphi({\bf
z})]
\end{equation}
that maps the zero sections $0_E\subset {\rm D}_1(E)$ onto
$Z_0\subset{\rm P}(E\oplus\C)$. Here ${\rm D}_1(E)$ is defined by
(\ref{e:1.3}), i.e., ${\rm D}_1(E)={\rm U}(E)\times_{U(n)}
B^{2n}(1)$.

\begin{theorem}\label{th:3.2}
For each given connection $1$-form $A$ on ${\mathcal F}:{\rm
U}(E)\to B$, there exist a small
$\varepsilon_0=\varepsilon_0({\mathcal F}, A)>0$  and two smooth
families of (compatible) symplectic forms
$\{\omega_{A\varepsilon}\,|\, 0<\varepsilon<\varepsilon_0\}$
\footnote{Hereafter the subscript $A\varepsilon$ does not mean to
multiply $A$ by $\varepsilon$!} on ${\rm D}_1(E)$, and
$\{\Omega_{A\varepsilon}\,|\, 0<\varepsilon<\varepsilon_0\}$ on
${\rm P}(E\oplus\C)$ such that
 the following properties hold:
\begin{itemize}
\item[($1^\circ$)] The maps $(B^{2n}(1),\varepsilon^2\omega_{\rm
std})\hookrightarrow ({\rm D}_1(E), \omega_{A\varepsilon})$ and
$(\CP^n,\varepsilon^2\omega_{\rm FS})\hookrightarrow ({\rm
P}(E\oplus\C), \Omega_{A\varepsilon})$
 are all symplectic
embeddings.

\item[($2^\circ$)] The map $\varphi_{\mathcal P}$ in
(\ref{e:3.10}) is a symplectic embedding from $({\rm D}_1(E),
\omega_{A\varepsilon})$ to $({\rm P}(E\oplus\C),
\Omega_{A\varepsilon})$.

 \item[($3^\circ$)] The zero section $0_E\subset E$ is a symplectic
 submanifold in $({\rm D}_1(E), \omega_{A\varepsilon})$ with symplectic
 normal bundle $(E,\varepsilon^2\bar\omega)\to B$. More precisely, for any $b\in
 B$,
\begin{equation}\label{e:3.10.1}
\omega_{A\varepsilon}|_{{\rm
D}_1(E)_b}=\varepsilon^2(\omega_{\rm std})_b\quad{\rm
and}\quad\Omega_{A\varepsilon}|_{{\rm
P}(E\oplus\C)_b}=\varepsilon^2(\omega_{\rm FS})_b,
\end{equation}
where $(\omega_{\rm std})_b$ is defined by (\ref{e:3.1}), and
$(\omega_{\rm FS})_b$ comes from ${\rm P}(E\oplus\C)$ being a
symplectic fibration with symplectic fibre $(\CP^n,\omega_{\rm
FS})$. So the symplectic forms $\omega_{A\varepsilon}$ and
$\Omega_{A\varepsilon}$ are compatible with the symplectic
fibrations.

 \item[($4^\circ$)]  The zero section $Z_0$ and
 the divisor $Z_\infty$ at infinity are
 symplectic submanifolds  in $({\rm
P}(E\oplus\C),\Omega_{A\varepsilon})$.

 \item[($5^\circ$)] If $F_0$ is a Hamiltonian $G$-subspace of
$(B^{2n}(1),\omega_{\rm std})$ (resp., $(\CP^n,\omega_{\rm FS})$),
then ${\rm U}(E)\times_{U(n)}F_0$ is a symplectic submanifold in
$({\rm D}_1(E), \omega_{A\varepsilon})$ (resp., $({\rm
P}(E\oplus\C), \Omega_{A\varepsilon})$).

\item[($6^\circ$)]  For any compact symplectic submanifold
$B^\circ\subset B$, let ${\rm U}(E)^\circ={\rm U}(E)|_{B^\circ}$
and $\{\omega_{A\varepsilon}^\circ\,|\,
0<\varepsilon<\varepsilon_0^\circ\}$ (resp.,
$\{\Omega_{A\varepsilon}^\circ\,|\,0<\varepsilon<\varepsilon_0^\circ
\}$) be the corresponding family of symplectic forms on ${\rm
D}_1(E)^\circ:={\rm U}(E)^\circ\times_{U(n)}B^{2n}(1)$ (resp.,
${\rm P}(E\oplus\C)^\circ:={\rm U}(E)^\circ\times_{U(n)}\CP^n$)
obtained by the restriction connection form $A|_{{\rm
U}(E)^\circ}$.
 Then  $\omega_{A\varepsilon}^\circ$ (resp.,
$\Omega_{A\varepsilon}^\circ$) is equal to the restriction of
 $\omega_{A\varepsilon}$ (resp., $\Omega_{A\varepsilon}$)
to ${\rm D}_1(E)^\circ:={\rm U}(E)^\circ\times_{U(n)}B^{2n}(1)$
(resp., ${\rm P}(E\oplus\C)^\circ:={\rm
U}(E)^\circ\times_{U(n)}\CP^n$) for each
$0<\varepsilon<\min\{\varepsilon_0^\circ, \varepsilon_0\}$.

\item[($7^\circ$)] If $P$ is  the trivial principal bundle
$P=B\times U(n)$ and $A$ is taken as  the connection form $A_{\rm
can}$ of the canonical flat connection in $P$ given by
Theorem~\ref{th:2.5}(iv),
  then for each $0<\varepsilon<\varepsilon_0(P, A_{\rm can})$,
 the symplectic form $\omega_{A_{{\rm can}\varepsilon}}$
(resp., $\Omega_{A_{{\rm can}\varepsilon}}$)
 on $P\times_{U(n)}B^{2n}(1)=B\times B^{2n}(1)$ (resp.,
 $P\times_{U(n)}\CP^n=B\times \CP^n$) is equal to
  $\beta\oplus (\varepsilon^2\omega_{\rm std})$ (resp.,
 $\beta\oplus (\varepsilon^2\omega_{\rm FS})$) for each
 $0<\varepsilon<\varepsilon_0(P, A_{\rm can})$.

 \item[($8^\circ$)] For each connection $1$-form $A$ on ${\rm
 U}(E)$,
let $\varepsilon_0=\varepsilon_0({\mathcal F}, A)>0$ be as above.
Then there exists a unique symplectic form  $\bar\omega_A$ on
${\rm
 D}_{\varepsilon_0}(E)$ such that for each
 $0<\varepsilon<\varepsilon_0$,
  \begin{equation}\label{e:3.10.2}
 \psi_\varepsilon: ({\rm
 D}_{1}(E), \omega_{A\varepsilon})\to ({\rm
 D}_{\varepsilon}(E), \bar\omega_{A}),\;(b,v)\mapsto (b,\varepsilon
 v)
 \end{equation}
 is a symplectomorphism.

 \item[($9^\circ$)]
For any two different connection forms $A_1$ and $A_2$ on
${\mathcal F}:{\rm U}(E)\to B$,
$0<\varepsilon\ll\min\{\varepsilon_0({\mathcal F},
A_1),\varepsilon_0({\mathcal F}, A_2)\}$, and $0<\delta\le 1$,
there exist fibre bundle isomorphisms
\begin{equation}\label{e:3.10.3}
\varphi^\delta_\varepsilon: {\rm D}_\delta(E)\to {\rm
D}_\delta(E)\quad{\rm and}\quad \Phi_\varepsilon: {\rm
P}(E\oplus\C)\to {\rm P}(E\oplus\C)
 \end{equation}
 such that $(\varphi^\delta_\varepsilon)^\ast\omega_{A_2\varepsilon}=
 \omega_{A_1\varepsilon}$ and
 $\Phi_\varepsilon^\ast\Omega_{A_2\varepsilon}=\Omega_{A_1\varepsilon}$.
 \end{itemize}
\end{theorem}

\begin{proof} The proofs can be obtained by
Theorem~\ref{th:2.5} directly. We only outline them. By the
explicit equivalence between ${\rm U}(E)\times_{U(n)}\C^{n}\to B$
and $E\to B$ given by
\begin{equation}\label{e:3.11}
\Xi: [p, (z_1,\cdots, z_n)]\mapsto z_1v_1+\cdots + z_nv_n,
\end{equation}
where $p=(v_1,\cdots,v_n)$ is a unitary frame of $E_b$, for any
$\epsilon>0$, we can write
\begin{equation}\label{e:3.11.1}
{\rm
D}_\epsilon(E)=\{(b,v)\in E\,|\, g_{J_E}(v,v)<\epsilon^2\}.
\end{equation}
 For a given connection $A$ on the principal bundle
${\mathcal F}:{\rm U}(E)\to B$, let ${\mathcal W}={\mathcal W}(A)$
be the largest open neighborhood of zero in ${\frak u}(n)^\ast$ so
that the corresponding minimal coupling form $\delta_A$ is
symplectic in ${\rm U}(E)\times{\mathcal W}$. Since $Cl(B^{2n}(1))$
and $\CP^n$ are compact, by (\ref{e:3.4}) and
Proposition~\ref{prop:3.1}, we can choose  \footnote{Here the choice
of $\varepsilon_0(A,{\mathcal F})$ shows that it is not canonical.}
$0<\varepsilon_0=\varepsilon_0(A,{\mathcal F})\le 1$ such that
\begin{equation}\label{e:3.12}
 \varepsilon^2\mu_{U(n)}(B^{2n}(1))\subset {\mathcal W}\quad{\rm
and}\quad
 \varepsilon^2\mu''_{U(n)}(\CP^{n})\subset {\mathcal
W}
\end{equation}
for all $\varepsilon\in (0,\varepsilon_0)$.
 (For example, ones
can take $\varepsilon_0$  to be the supremum of $\varepsilon>0$
satisfying the  inclusion relations in (\ref{e:3.12}).) Now
applying Theorem~\ref{th:2.4} to the cases $P={\rm U}(E)$ and
$(F,\sigma)=(B^{2n}(1), \varepsilon^2\omega_{\rm std})\;{\rm or}\;
(\CP^n,\varepsilon^2\omega_{\rm FS})$ we immediately get
($1^\circ$).

Note that (\ref{e:3.9}) means that the symplectic embedding
$\varphi$ in (\ref{e:3.8}) commutes with the Hamiltonian actions
in (\ref{e:1.1}) and (\ref{e:3.5}). Condition ($2^\circ$) follows
from Theorem~\ref{th:2.5}(i).

The conclusion in ($3^\circ$) is easily derived from
Theorem~\ref{th:2.5}(ii). It precisely says\linebreak
$(T_{(b,0)}0_E)^{\omega_{A\varepsilon}}=(E_b,\varepsilon^2\bar\omega_b)$
for any $b\in B$.

To see ($4^\circ$), note that $\CP^{n-1}=\{[z_0,\cdots,
z_{n-1},0]\in\CP^n\}\subset\CP^n$ is a symplectic submanifold in
$(\CP^n,\omega_{\rm FS})$ that is invariant under the action in
(\ref{e:3.5}). By Theorem~\ref{th:2.5}(i), $Z_\infty$ is a
symplectic submanifold in $({\rm
P}(E\oplus\C),\Omega_{A\varepsilon})$. Since the symplectic
embedding $\varphi_{\mathcal P}$ in ($2^\circ$) maps the zero
section $0_E$ onto $Z_0={\rm P}(\{0\}\oplus\C)$, $Z_0$ is a
symplectic submanifold in $({\rm
P}(E\oplus\C),\Omega_{A\varepsilon})$.

Theorem~\ref{th:3.2}($5^\circ$) is a direct consequence of the
final conclusion of Theorem~\ref{th:2.5}(i).
Theorem~\ref{th:3.2}($6^\circ$), ($7^\circ$) follows from
Theorem~\ref{th:2.5}(iii), (iv), respectively. Condition
($8^\circ$) can easily be obtained by Theorem~\ref{th:2.5}(i).
Finally, we prove ($9^\circ$).  Let ${\mathcal
W}_1^\star\subset{\mathcal W}(A_1)$ and ${\mathcal
W}_2^\star\subset{\mathcal W}(A_2)$ be open neighborhoods of zero
in ${\frak u}(n)^\ast$ such that $({\rm U}(E)\times{\mathcal
W}_1^\star, \delta_{A_1})$ is symplectomorphic to $({\rm
U}(E)\times{\mathcal W}_2^\star, \delta_{A_2})$. Then for
$0<\varepsilon\ll\min\{\varepsilon_0({\mathcal F},
A_1),\varepsilon_0({\mathcal F}, A_2)\}$,  it holds that
$\varepsilon^2\mu_{U(n)}(Cl(B^{2n}(1)))\subset {\mathcal W}^\star$
and $\varepsilon^2\mu''_{U(n)}(\CP^{n})\subset {\mathcal
W}^\star$. Under this case Theorem~\ref{th:2.5}(v) gives the
desired results directly. Theorem~\ref{th:3.2} are proved.
\end{proof}

 Since $\omega_{A\varepsilon}$ is compatible with symplectic
 fibration, by Theorem~\ref{th:3.2}($3^\circ$),  we can  write
\begin{equation}\label{e:3.13}
\left.\begin{array}{ll}
\omega_{A\varepsilon}=\pi^\ast\beta+\tau_{A\varepsilon},\quad
\tau_{A\varepsilon}|_{0_E}=0,\\
 \tau_{A\varepsilon}|_{{\rm D}_1(E)_b}=\varepsilon^2(\omega_{\rm
std})_b\quad\forall b\in B.
\end{array}\right.
\end{equation}
Similarly, by Theorem~\ref{th:3.2}($4^\circ$), we easily get that
$\Omega_{A\varepsilon}|_{Z_0}=\beta$ after identifying $Z_0\equiv
0_E\equiv B$. It follows that we can also write
\begin{equation}\label{e:3.14}
\left.\begin{array}{ll} \Omega_{A\varepsilon}={\mathcal
P}^\ast\beta+\Gamma_{A\varepsilon},\quad
\Gamma_{A\varepsilon}|_{Z_0}=0,\\
 \Gamma_{A\varepsilon}|_{{\rm
P}(E\oplus\C)_b}=\varepsilon^2(\omega_{\rm FS})_b\quad\forall b\in
B.
\end{array}\right.
\end{equation}
Note that the almost complex structure ${\bf J}$ on ${\rm
P}(E\oplus\C)$ constructed by (\ref{e:4.6}) is not necessarily
$\Omega_{A\varepsilon}$-tamed. However, we can show that for
sufficiently small $t>0$, the closed 2-form ${\mathcal P}^\ast\beta+
t\Gamma_{A\varepsilon}$ is also symplectic and tame this ${\bf J}$
(see Lemma~\ref{lem:4.2}). Hence we are led to the following
strengthened version of Theorem~\ref{th:3.2}, whose precise
statement is needed in the proof of Theorem~\ref{th:4.5}.

\begin{theorem}\label{th:3.3}
Under the assumption of Theorem~\ref{th:3.2}, for each
$0<\varepsilon<\varepsilon_0$,  there exists a small $t_0=t_0(A,
\varepsilon)>0$ such that for each $0<t<t_0$, the form
\begin{equation}\label{e:3.15}
\omega^t_{A\varepsilon}:=\pi^\ast\beta+
t\tau_{A\varepsilon}\quad({\rm resp.,}\quad
\Omega^t_{A\varepsilon}:={\mathcal P}^\ast\beta+
t\Gamma_{A\varepsilon})
\end{equation}
is a symplectic form on  ${\rm D}_1(E)$ (resp. ${\rm
P}(E\oplus\C)$), where $\tau_{A\varepsilon}$ and
$\Gamma_{A\varepsilon}$ are given by (\ref{e:3.13}) and
(\ref{e:3.14}), respectively. Moreover, they also satisfy the
following properties.
\begin{itemize}
\item[(i)] The map $(B^{2n}(1), t\varepsilon^2\omega_{\rm
std})\hookrightarrow ({\rm D}_1(E), \omega^t_{A\varepsilon})$ and
$(\CP^n, t\varepsilon^2\omega_{\rm FS})\hookrightarrow ({\rm
P}(E\oplus\C),\Omega^t_{A\varepsilon})$ are symplectic embeddings.

\item[(ii)] The map $\varphi_{\mathcal P}$ in (\ref{e:3.10})  is a
symplectic embedding from $({\rm D}_1(E),
\omega^t_{A\varepsilon})$ into $({\rm P}(E\oplus\C),
\Omega^t_{A\varepsilon})$.

 \item[(iii)] The zero section $0_E\subset E$ is a symplectic
 submanifold in $({\rm D}_1(E), \omega^t_{A\varepsilon})$ with symplectic
 normal bundle $(E,t\varepsilon^2\bar\omega)\to B$.

 \item[(iv)]  The zero section $Z_0={\rm P}(\{0\}\oplus \C)$ is
 a symplectic submanifold  in $({\rm
P}(E\oplus\C), \Omega^t_{A\varepsilon})$.

\item[(v)] Let the closed two-form
  $\bar\omega_A^t$ be defined by (\ref{e:1.5}).
Then for any $0<\varepsilon<\varepsilon_0(A, {\mathcal F})$ and
$0<t<t_0(A,\varepsilon)$,  $\bar\omega_A^t$ is symplectic in ${\rm
D}_{\varepsilon}(E)$ and
\begin{equation}\label{e:3.15.1}
 \psi_\varepsilon: ({\rm
 D}_{1}(E), \omega^t_{A\varepsilon})\to ({\rm
 D}_{\varepsilon}(E), \bar\omega^t_{A}),\;(b,v)\mapsto (b,\varepsilon
 v)
 \end{equation}
 is also a symplectomorphism.

\item[(vi)]  For any compact symplectic submanifold
$B^\circ\subset B$ of codimension zero, as in
Theorem~\ref{th:3.2}($6^\circ$), let $\omega_{A\varepsilon}^\circ$
(resp., $\Omega_{A\varepsilon}^\circ$) be the symplectic form on
${\rm D}_1(E)^\circ$ (resp., ${\rm P}(E\oplus\C)^\circ$) for
$0<\varepsilon<\varepsilon_0^\circ$. As in (\ref{e:3.15}), there
exist a small $t^\circ_0(A,\varepsilon)>0$ and two families of
symplectic forms,
\begin{equation}\label{e:3.15.2} \omega^{\circ
t}_{A\varepsilon}=(\pi|_{{\rm D}_1(E)^\circ}) ^\ast\beta +
t\tau^\circ_{A\varepsilon}\quad{\rm and}\quad \Omega^{\circ
t}_{A\varepsilon}=({\mathcal P}|_{{\rm
P}(E\oplus\C)^\circ})^\ast\beta + t\Gamma^\circ_{A\varepsilon}
\end{equation}
for $0<t<t_0^\circ(A,\varepsilon)$, where
$\omega^{\circ}_{A\varepsilon}=(\pi|_{{\rm D}_1(E)^\circ})
^\ast\beta + \tau^\circ_{A\varepsilon}$ and $\Omega^{\circ
}_{A\varepsilon}=({\mathcal P}|_{{\rm
P}(E\oplus\C)^\circ})^\ast\beta \\ + \Gamma^\circ_{A\varepsilon}$.
Then for each $0<\varepsilon<\min\{\varepsilon_0^\circ,
\varepsilon_0\}$ and $0<t<\min\{t_0, t_0^\circ\}$, the symplectic
form $\omega_{A\varepsilon}^{\circ t}$ (resp.,
$\Omega_{A\varepsilon}^{\circ t}$) is equal to the restriction of
 $\omega^t_{A\varepsilon}$ (resp., $\Omega^t_{A\varepsilon}$)
to ${\rm D}_1(E)^\circ$ (resp., ${\rm P}(E\oplus\C)^\circ$).

 \item[(vii)] Under the assumptions of Theorem~\ref{th:3.2}($9^\circ$),
for any \linebreak $0<t<\min\{t_0(A_1,\varepsilon),
t_0(A_2,\varepsilon)\}$ and $\delta\in (0, 1]$,
\begin{equation}\label{e:3.15.3}
\left.\begin{array}{ll}
& \varphi^\delta_\varepsilon: ({\rm
 D}_\delta(E), \omega^t_{A_1\varepsilon})\to ({\rm
 D}_\delta(E), \omega^t_{A_2\varepsilon}),\\
  &\Phi_\varepsilon:({\rm P}(E\oplus\C),
\Omega^t_{A_1\varepsilon})\to({\rm P}(E\oplus\C),
\Omega^t_{A_2\varepsilon})\end{array}\right.
\end{equation}
are all symplectomorphisms.

\item[(viii)] If $P=B\times U(n)$ and $A=A_{\rm can}$ are as in
Theorem~\ref{th:3.2}($7^\circ$), then for any
$0<\varepsilon<\varepsilon_0(P, A_{\rm can})$, one can take
$t_0(A_{\rm can},\varepsilon)=\infty$ and $\bar\omega^t_{A_{{\rm
can}\varepsilon}}=\beta\oplus(t\varepsilon^2\omega_{\rm std})$ for
each $t>0$.
\end{itemize}
\end{theorem}

\begin{proof}We firstly prove that
$\Omega^t_{A\varepsilon}={\mathcal P}^\ast\beta+
t\Gamma_{A\varepsilon}$ are symplectic forms on ${\rm
P}(E\oplus\C)$ for sufficiently small $t>0$. For any $x\in{\rm
P}(E\oplus\C)$ let
\begin{equation}\label{e:3.15.4} {\mathcal
V}_x:={\rm Ker}d{\mathcal P}(x)\subset T_x{\rm P}(E\oplus\C).
\end{equation}
Then it is equal to $T_x{\rm P}(E\oplus\C)_b$, where $b={\mathcal
P}(x)$. Since each fibre is a symplectic submanifold, the subspace
of $ T_x{\rm P}(E\oplus\C)$,
\begin{eqnarray}\label{e:3.15.5}
{\mathcal H}_x:=({\mathcal V}_x)^{\Omega_{A\varepsilon}}=\{X\in
T_x{\rm P}(E\oplus\C)\,|\, \Omega_{A\varepsilon}(X,Y)=0\;\forall
Y\in{\mathcal
V}_x\}\nonumber\\
=\{X\in T_x{\rm P}(E\oplus\C)\,|\,
\Gamma_{A\varepsilon}(X,Y)=0\;\forall Y\in{\mathcal V}_x\}
\end{eqnarray}
is not only symplectic, but also a horizontal complement of
${\mathcal V}_x$, that is,
\begin{equation}\label{e:3.16}
T_x{\rm P}(E\oplus\C)={\mathcal H}_x\oplus{\mathcal V}_x
\end{equation}
is a symplectic direct sum decomposition and the projection
\begin{equation}\label{e:3.17}
d{\mathcal P}(x):{\mathcal H}_x\to T_bB
\end{equation}
 is a bijection. By Lemma~\ref{lem:4.2}, we
 have an almost complex structure ${\bf J}$ on
${\rm P}(E\oplus\C)$ and a small $\bar t>0$ such that
$\Omega_{A\varepsilon}^t(X, {\bf J}X)>0$ for any nonzero $X\in
T{\rm P}(E\oplus\C)$ and $t\in (0, \bar t]$. So, such
$\Omega_{A\varepsilon}^t$ must be nondegenerate. Then the desired
$t_0$ can be taken as the supremum of $\bar t>0$ for which all
$\Omega_{A\varepsilon}^t$ are nondegenerate for any $t\in (0, \bar
t]$. The conclusions for $\omega^t_{A\varepsilon}$ can be proved
 in the same way (by constructing an almost complex structure
 on the compact bundle $Cl({\rm D}_1(E))$ as in (\ref{e:4.6})).
 Of course we can shrink $t_0$ if necessary.

Conditions (i), (iii), (iv) and (viii) are obvious. To see (ii),
note that $\varphi^\ast_{\mathcal P}({\mathcal
P}^\ast\beta)=\pi^\ast(\varphi^\ast_{\mathcal P}\beta)$. So
$\varphi_{\mathcal
P}^\ast\Omega_{A\varepsilon}=\omega_{A\varepsilon}$ if and only if
$\varphi_{\mathcal
P}^\ast\Gamma_{A\varepsilon}=\tau_{A\varepsilon}$. The desired
conclusion follows immediately.

 To prove (v), we still use $\pi$ to denote the bundle
projections from ${\rm D}_{1}(E)$ and ${\rm D}_{\varepsilon}(E)$
to $B$. Then $\pi\circ\psi_\varepsilon=\pi$, and hence $d\pi\circ
d\psi_\varepsilon=d\pi$. In particular, we have
\begin{equation}\label{e:3.18}
\psi_\varepsilon^\ast(\pi^\ast\beta)=\pi^\ast\beta
\end{equation}
 For any
$x\in{\rm D}_{1}(E)$ and $y\in {\rm D}_{\varepsilon}(E)$ let
\begin{equation}\label{e:3.18.1}
\left.\begin{array}{ll}
&{\bf V}_x:={\rm Ker}(d\pi(x)),\quad {\bf
H}_x:=({\bf
V}_x)^{\omega_{A\varepsilon}},\\
&{\bf V}'_y:={\rm Ker}(d\pi(y)),\quad {\bf H}'_y:=({\bf
V}'_y)^{\bar\omega_A}.\end{array}\right.
\end{equation}
Then $T_x{\rm D}_{1}(E)={\bf H}_x\oplus{\bf V}_x$ and $T_y{\rm
D}_{\varepsilon}(E)={\bf H}'_y\oplus{\bf V}'_y$. Clearly,
$d\psi_\varepsilon(x)({\bf V}_x)={\bf V}'_{\psi_\varepsilon(x)}$.
Since $\psi_\varepsilon^\ast\bar\omega_A=\omega_{A\varepsilon}$,
one easily derives that $d\psi_\varepsilon(x)({\bf H}_x)={\bf
H}'_{\psi_\varepsilon(x)}$. These imply that
\begin{equation}\label{e:3.19}
\psi_\varepsilon^\ast(\bar\omega_A-\pi^\ast\beta)=\omega_{A\varepsilon}-\pi^\ast\beta.
\end{equation}
Hence (\ref{e:3.18}) and (\ref{e:3.19}) together show that
$\psi_\varepsilon^\ast\bar\omega_A^t=\omega^t_{A\varepsilon}$ for
any $t$. But $\omega^t_{A\varepsilon}$ is symplectic in ${\rm
D}_1(E)$ for any $t\in (0, t_0)$. Hence $\bar\omega_A^t$ is
symplectic in ${\rm D}_{\varepsilon}(E)$ for each
$0<t<t_0(A,\varepsilon)$. The desired conclusion is proved.

Condition (vi) easily follows from Theorem~\ref{th:3.2}($6^\circ$)
and the proof of (v) above.

Finally, since $\pi\circ\varphi_\varepsilon=\pi$ and ${\mathcal
P}\circ\Phi_\varepsilon={\mathcal P}$, as in the proof of (v) ones
can
 derive (vii)  from Theorem~\ref{th:3.2}($9^\circ$) easily.\end{proof}

\section{Proof of Main Result}
\setcounter{equation}{0}

We briefly review the definition of the pseudo symplectic capacity
of the Hofer-Zehnder type introduced by the author in \cite{Lu2}.
See [17, 19, 20, 18] for more estimates and applications.

For a connected symplectic manifold $(M,\omega)$ of dimension at
least $4$ and two nonzero homology classes
$\alpha_0,\alpha_\infty\in H_\ast(M; \Q)$,  we say that a smooth
function  $H:M\to\R$ is
\emph{$(\alpha_0,\alpha_\infty)$-admissible}
  (resp.,\  $(\alpha_0,\alpha_\infty)^\circ$-\emph{admissible}) if
there exist two compact submanifolds $P$ and $Q$ of $M$ with
connected smooth boundaries and of codimension zero such that the
following condition groups $(1), (2), (3), (4), (5)$, and $(6)$
(resp.,\ $(1), (2), (3), (4), (5)$, and $(6^\circ))$ hold:
\begin{itemize}
\item[(1)] $P\subset \Int (Q)$ and $Q\subset \Int (M)$;

\item[(2)]$H|_P=0$ and $H|_{M \setminus \Int (Q)}=\max H$;

\item[(3)] $0\le H\le \max H$; \item[(4)] There exist chain
representatives of $\alpha_0$ and $\alpha_\infty$,
            still denoted by $\alpha_0$ and $\alpha_\infty$, such that
            ${\it supp}(\alpha_0)\subset \Int (P)$ and
            ${\it supp}(\alpha_\infty)\subset M\setminus Q$;

\item[(5)] There are no critical values in $(0,\varepsilon)\cup
(\max H-\varepsilon, \max H)$
    for a small $\varepsilon=\varepsilon(H)>0$;

  \item[(6)] The
Hamiltonian system $\dot x=X_H(x)$ on $M$ has no nonconstant fast
periodic solutions; \item[(6$^\circ$)] The Hamiltonian system
$\dot x=X_H(x)$ on $M$ has no nonconstant contractible fast
periodic solutions.
\end{itemize}

Let ${\mathcal H}_{ad}(M,\omega;\alpha_0,\alpha_\infty)$ and
${\mathcal H}_{ad}^\circ(M,\omega;\alpha_0,\alpha_\infty)$ denote
the sets of
 $(\alpha_0,\alpha_\infty)$-admissible functions and
$(\alpha_0,\alpha_\infty)^\circ$-admissible ones, respectively. In
\cite{Lu2}, we defined
\begin{equation}\label{e:4-2}
\begin{array}{ll}
C_{HZ}^{(2)}(M,\omega;\alpha_0,\alpha_\infty):=\sup\left\{\max
H\,|
\, H\in {\mathcal H}_{ad}(M,\omega;\alpha_0,\alpha_\infty)\right\},\\
C_{HZ}^{(2\circ)}(M,\omega;\alpha_0,\alpha_\infty):=\sup\left\{\max
H\,| \, H\in {\mathcal
H}_{ad}^\circ(M,\omega;\alpha_0,\alpha_\infty)\right\}.
\end{array}
\end{equation}
They were, respectively, called the \emph{pseudo symplectic
capacity of the Hofer-Zehnder type} and the
$\pi_1$-\emph{sensitive pseudo symplectic capacity of the
Hofer-Zehnder type}. In particular, we get a genuine symplectic
capacity
\begin{equation}\label{e:4-1}
C_{HZ}^{(2)}(M,\omega):=C_{HZ}^{(2)}(M,\omega; pt, pt)
\end{equation}
and a $\pi_1$-sensitive symplectic capacity
\begin{equation}\label{e:4.0}
C_{HZ}^{(2\circ)}(M,\omega):=C_{HZ}^{(2\circ)}(M,\omega; pt, pt).
\end{equation}
We also showed in \cite[Lemma 1.4]{Lu2} that there exist the
following relations among them, the usual Hofer-Zehnder capacity
$c_{HZ}$ and the $\pi_1$-sensitive Hofer-Zehnder capacity
$c^\circ_{HZ}$:
\begin{equation}\label{e:4.1}
C^{(2)}_{HZ}(M, \omega)=c_{HZ}(M, \omega)\quad{\rm and}\quad
C_{HZ}^{(2\circ)}(M, \omega)=c^\circ_{HZ}(M,\omega)
\end{equation}
\emph{if a symplectic manifold $(M,\omega)$ is either closed or
satisfies the condition that for each compact subset $K\subset
M\setminus\partial M$, there exists a compact submanifold $W\subset
M$ with connected boundary and of codimension zero such that
$K\subset W$.}

 For a closed symplectic manifold $(M,\omega)$ and
$\alpha_0,\alpha_\infty\in H_\ast(M;\Q)$, let
\begin{equation}\label{e:4.1.1}
{\rm
GW}_g(M,\omega;\alpha_0,\alpha_\infty)\in (0, +\infty]
\end{equation}
be the infimum of the $\omega$-areas $\omega(A)$ of the homology
classes $A\in H_2(M;\Z)$ for which the Gromov-Witten invariant
$\Psi_{A, g,
m+2}(C;\alpha_0,\alpha_\infty,\beta_1,\cdots,\beta_m)\ne 0$ for
some homology classes $\beta_1,\cdots,\beta_m\in H_\ast(M;\Q)$ and
$C\in  H_\ast(\overline{\mathcal M}_{g, m+2};\Q)$ and an integer
$m \ge 1$. Here for a given class $A\in H_2(M;\Z)$, the
Gromov-Witten invariant of genus $g$ and with $m+2$ marked points
is a homomorphism
\begin{equation}\label{e:4.1.2.0}
\Psi_{A, g, m+2}: H_\ast(\overline{\mathcal M}_{g, m+2};\Q)\times
H_\ast(M;\Q)^{m+2}\to \Q,
\end{equation}
 the reader may refer to
\cite{LiuT2, Lu6} for details. (In the latter paper, we used the
cohomology and denoted by ${\mathcal G}{\mathcal W}$ the
GW-invariants. It is easily translated into the homology while $M$
is a closed manifold.)  We also define
\begin{equation}\label{e:4.1.2} {\rm
GW}(M,\omega;\alpha_0,\alpha_\infty):= \inf \left\{ {\rm
GW}_g(M,\omega;\alpha_0,\alpha_\infty)\,|\, g\ge 0 \right\}\in [0,
+\infty].
\end{equation}
Based on \cite{LiuT2}, we proved in \cite[Theorem 1.10]{Lu2} that
\begin{eqnarray}
&&C_{HZ}^{(2)}(M,\omega;\alpha_0,\alpha_\infty)\le {\rm
GW}(M,\omega;\alpha_0,\alpha_\infty),
\label{e:4.2}\\
&&C_{HZ}^{(2\circ)}(M,\omega;\alpha_0,\alpha_\infty)\le {\rm
GW}_0(M,\omega;\alpha_0,\alpha_\infty)\label{e:4.3}
\end{eqnarray}
for any closed symplectic manifold $(M,\omega)$ of dimension $\dim
M \ge 4$ and homology classes $\alpha_0,\alpha_\infty\in
H_\ast(M;\Q)\setminus\{0\}$. The following proposition lists
partial results in \cite[Proposition 1.7]{Lu2}, which are needed
in the following arguments.

\begin{proposition}  \label{prop:4.1}
  Let $W\subset \Int (M)$ be a smooth compact submanifold of codimension
  zero and with connected boundary such that
the homology classes $\alpha_0,\alpha_\infty\in
H_\ast(M;\Q)\setminus\{0\}$ have representatives supported in
$\Int (W)$ and $\Int (M)\setminus W$, respectively. Denote by
$\tilde\alpha_0\in H_\ast(W; \Q)$ and $\tilde\alpha_\infty\in
H_\ast(M\setminus W; \Q)$ the nonzero homology classes determined
by them. Then
\begin{equation}  \label{e:4.3.1}
C_{HZ}^{(2)}(W,\omega;\tilde\alpha_0, pt) \le
C_{HZ}^{(2)}(M,\omega;\alpha_0,\alpha_\infty)
\end{equation}
 and,
in particular, one has
\begin{equation}  \label{e:4.4}
c_{HZ}(W,\omega)=C^{(2)}_{HZ}(W,\omega)\le C_{HZ}^{(2)}(M,\omega;
pt, \alpha)
\end{equation}
for any $\alpha\in H_\ast(M; \Q)\setminus\{0\}$ with
representative supported in $\Int (M) \setminus W$. If the
inclusion $W\hookrightarrow M$ induces an injective homomorphism
$\pi_1(W)\to\pi_1(M)$, then
\begin{equation}  \label{e:4.4.1}
C_{HZ}^{(2\circ)}(W,\omega;\tilde\alpha_0, pt) \le
C_{HZ}^{(2\circ)}(M,\omega;\alpha_0,\alpha_\infty)
\end{equation}
 and corresponding to (\ref{e:4.4}), one has
\begin{equation}  \label{e:4.5}
c_{HZ}^\circ(W,\omega)=C_{HZ}^{(2\circ)}(W,\omega)\le
C_{HZ}^{(2\circ)}(M,\omega; pt, \alpha).
\end{equation}
\end{proposition}

We now construct a class of almost complex structures on ${\rm
P}(E\oplus\C)$. The readers will see why we need
Theorem~\ref{th:3.3}. Note that every given Riemannian metric
${\rm P}(E\oplus\C)$ naturally restricts to a natural Riemannian
metric $g_b$ on the fibre ${\rm P}(E\oplus\C)_b$ for each $b\in
B$. Using the standard method, (cf. [29, page 64]) one gets a
compatible almost complex structure $J_b\in {\mathcal J}({\rm
P}(E\oplus\C)_b,(\omega_{\rm FS})_b)$. Clearly, $J_b$ smoothly
depends on $b\in B$. Take another compatible almost complex
structure $J_B\in {\mathcal J}(B,\beta)$. By (\ref{e:3.17}), we
can obtain its horizontal lift $\widetilde J_B$ to ${\mathcal H}$:
\begin{equation}  \label{e:4.5.1}
(\widetilde J_B)_x:{\mathcal H}_x\to {\mathcal H}_x,\;X\mapsto
(d{\mathcal P}(x))^{-1}\circ(J_B)_b\circ d{\mathcal P}(x)(X),
\end{equation}
where $b={\mathcal P}(x)$. Then using the decomposition
(\ref{e:3.16}), we may define an almost complex structure ${\bf
J}$ on ${\rm P}(E\oplus\C)$ as follows:
\begin{eqnarray}\label{e:4.6}
&&{\bf J}_x:T_x{\rm P}(E\oplus\C)={\mathcal H}_x\oplus{\mathcal
V}_x\to T_x{\rm
P}(E\oplus\C),\\
&&\qquad X^h\oplus X^v\mapsto(\widetilde J_B)_x(X^h)\oplus
(J_b)_x(X^v).\nonumber
\end{eqnarray}

\begin{lemma}\label{lem:4.2}
${\bf J}$ is $\Omega^t_{A\varepsilon}$-tamed for sufficiently
small $t>0$.
\end{lemma}

\begin{proof} For any nonzero $X=X^h+ X^v\in
T_x{\rm P}(E\oplus\C)={\mathcal H}_x\oplus{\mathcal V}_x$, the
direct computation yields
\begin{equation}\label{e:4.6.1}
\left.\begin{array}{ll}
 \Omega^t_{A\varepsilon}(X, {\bf
J}X)\!\!\!\!\!\!&=\Omega^t_{A\varepsilon}\bigl(X^h+X^v,
(\widetilde J_B)_xX^h+
(J_b)_xX^v\bigr)\\
&={\mathcal P}^\ast\beta\bigl(X^h+X^v, (\widetilde J_B)_xX^h+
(J_b)_xX^v\bigr)\\
&\quad + t\Gamma_{A\varepsilon}\bigl(X^h+ X^v,(\widetilde
J_B)_xX^h+
(J_b)_xX^v\bigr)\\
&=\beta\bigl(d{\mathcal P}(x)(X^h), d{\mathcal P}(x)((\widetilde
J_B)_xX^h)
\bigr)\\
&\quad + t\Gamma_{A\varepsilon}\bigl(X^h,(\widetilde
J_B)_xX^h\bigr)+t\varepsilon(\omega_{\rm FS})_b\bigl(X^v, (J_b)_xX^v\bigr)\\
&\quad + t\Gamma_{A\varepsilon}\bigl(X^h,(J_b)_xX^v\bigr) +
t\Gamma_{A\varepsilon}\bigl(X^v,(\widetilde
J_B)_xX^h\bigr)\\
&=\beta\bigl(d{\mathcal P}(x)(X^h), (J_B)_bd{\mathcal P}(x)X^h
\bigr)\\
&\quad + t\Gamma_{A\varepsilon}\bigl(X^h,(\widetilde
J_B)_xX^h\bigr)+ t\varepsilon(\omega_{\rm FS})_b\bigl(X^v,
(J_b)_xX^v\bigr).\end{array}\right.
\end{equation}
Here $\Gamma_{A\varepsilon}\bigl(X^h,(J_b)_xX^v\bigr)=0$ and
$\Gamma_{A\varepsilon}\bigl(X^v,(\widetilde J_B)_xX^h\bigr)=0$ are
because of $(J_b)_xX^v\in{\mathcal V}_x$, $(\widetilde
J_B)_xX^h\in{\mathcal H}_x$ and the equality above (\ref{e:3.16}),
that is, ${\mathcal H}_x=\{X\in T_x{\rm P}(E\oplus\C)\,|\,
\Gamma_{A\varepsilon}(X,Y)=0\;\forall Y\in{\mathcal V}_x\}$. Since
${\rm P}(E\oplus\C)$ is compact,  and by (\ref{e:3.17}) the
projection $d{\mathcal P}(x)$ is a bijection,  it easily follows
that for a given $\epsilon\in (0, 1)$, there exists a $0<\bar t\le
t_0$ such that
\begin{equation}\label{e:4.6.2} \left.\begin{array}{ll}
\beta\bigl(d{\mathcal P}(x)(X^h), (J_B)_bd{\mathcal P}(x)X^h\bigr)
+t\Gamma_{A\varepsilon}\bigl(X^h,(\widetilde J_B)_xX^h\bigr)\\
\qquad\ge (1-\epsilon)\beta\bigl(d{\mathcal P}(x)(X^h),
(J_B)_bd{\mathcal P}(x)X^h\bigr)\end{array}\right.
\end{equation}
for all $t\in (0, \bar t]$ and any $X\in T{\rm P}(E\oplus\C)$.
Note that $J_b$ is compatible with $(\omega_{\rm FS})_b$. The
desired claim is proved. \end{proof}

Note that the almost complex structure ${\bf J}$ in (\ref{e:4.6})
satisfies $d{\mathcal P}\circ{\bf J}=J_B\circ d{\mathcal P}$. That
is, the projection ${\mathcal P}$ is holomorphic with respect to the
almost complex structures ${\bf J}$ and $J_B$. In terms of [27,
Definition 2.8], we say that the almost complex structure ${\bf J}$
on ${\rm P}(E\oplus\C)$ is \emph{fibred} with respect to $J_B$. For
such an almost complex structure ${\bf J}$ on ${\rm P}(E\oplus\C)$,
any ${\bf J}$-holomorphic curve in ${\rm P}(E\oplus\C)$ representing
a homology class of fiber  must sit entirely in a fiber. By the
following remark, each $J_b$ can be chosen as a canonical one such
that for the corresponding almost complex structure ${\bf J}$, the
projection ${\mathcal P}:({\rm P}(E\oplus\C), {\bf J})\to (B, J_B)$
is an \emph{almost complex fibration with fibre} $(\CP^n, J_{\rm
FS})$ in the sense of [1, Definition 6.3.A].

\begin{remark}\label{rm:4.3}{\rm
The above  almost complex structure $J_b\in {\mathcal J}({\rm
P}(E\oplus\C)_b,(\omega_{\rm FS})_b)$ can be chosen as canonical one
as done in \cite{Bi}. Namely, a Hermitian vector bundle may yield a
canonical (almost) complex structure on each fibre of its projective
bundle.   For every $J_E\in{\mathcal J}(E,\bar\omega)$, one has a
corresponding almost complex structure on ${\rm P}(E)$, $J_{\rm
can}$ which restricts to  an almost complex structure $J_{\rm
can,b}$ on ${\rm P}(E)_b={\rm P}(E_b)$ for each $b\in B$.}
\end{remark}

Recall that in \cite[Definition 1.14]{Lu2}, a closed symplectic
manifold $(M,\omega)$ is called $g$-\emph{symplectic uniruled} if
the Gromov-Witten invariants
\begin{equation}\label{e:4.6.3}
\Psi_{A, g, m+2}(C; pt,\alpha,\beta_1,\cdots,\beta_m)\ne 0
\end{equation}
for some homology classes $A\in H_2(M;\Z)$, $\alpha,
\beta_1,\cdots, \beta_m\in H_\ast(M;\Q)$ and \linebreak $C\in
H_\ast(\overline{\mathcal M}_{g, m+2};\Q)$ and an integer $m \ge
1$. In particular, if $C$ can be chosen as the class of a point,
$(M,\omega)$ is said to be {\it strong} $g$-\emph{symplectic
uniruled}.

\begin{theorem}\label{th:4.4}
Under the assumptions of Theorem~\ref{th:3.3}, for the symplectic
forms  $\Omega^t_{A\varepsilon}$ in (\ref{e:3.15}), it holds that
every  symplectic manifold $({\rm
P}(E\oplus\C),\Omega^t_{A\varepsilon})$ is  strong $0$-symplectic
uniruled. More precisely, the GW-invariant
\begin{equation}\label{e:4.7}
\Psi_{L, 0, 3}(pt; pt, Z_0, Z_\infty)=1,
\end{equation}
where $L$ denotes the class of the line in the fiber of ${\rm
P}(E\oplus\C)$.
\end{theorem}

We put off its proof to the end of this paper.

\begin{theorem}\label{th:4.5}
 Under the assumptions of Theorem~\ref{th:4.4}, the pseudo symplectic capacity
\begin{equation}\label{e:4.8}
C_{HZ}^{(2\circ)}({\rm P}(E\oplus\C),\Omega^t_{A\varepsilon}; pt,
Z_0)\le\pi t\varepsilon^2
\end{equation}
for each $0<\varepsilon<\varepsilon_0(A,{\mathcal F})$ and $0<t<
t_0(A,\varepsilon)$. In particular, the $\pi_1$-sensitive
Hofer-Zehnder capacity
\begin{equation}\label{e:4.9}
c_{HZ}^\circ({\rm D}_{\varepsilon}(E), \bar\omega^t_A)=
c_{HZ}^\circ({\rm D}_{1}(E), \omega^t_{A\varepsilon})\le \pi
t\varepsilon^2
\end{equation}
for each $0<\varepsilon<\varepsilon_0({\mathcal F}, A)$  and
$0<t<t_0(A,\varepsilon)$.
 Furthermore, for a given connection form  $A'$ on
${\rm U}(E)$ which is {\rm flat} near some symplectically embedded
ball in $(B,\beta)$ of radius being the symplectic radius of
$(B,\beta)$, and each
\begin{equation}\label{e:4.9.1}
0<\varepsilon\ll \min\{\varepsilon_0({\mathcal F}, A),
\varepsilon_0({\mathcal F},
 A'), \varepsilon_0(P, A_{\rm can}), \sqrt{\frac{{\mathcal
 W}_G(B,\beta)}{\pi}}\},\end{equation}
 it holds that
\begin{equation}\label{e:4.10}
c({\rm D}_{\varepsilon}(E), \bar\omega^t_A)= c({\rm D}_{1}(E),
\omega^t_{A\varepsilon})=\pi t\varepsilon^2
\end{equation}
for each $0<t<\min\{t_0(A',\varepsilon), t_0(A,\varepsilon)\}$ and
$c={\mathcal W}_G$, $c_{HZ}$, and $c_{HZ}^\circ$.
 \end{theorem}

\begin{proof} {\bf Step 1.} {\it Proving
(\ref{e:4.9})}. Firstly, note that (\ref{e:4.3}) and (\ref{e:4.7})
directly yields (\ref{e:4.8}) because
$\Omega^t_{A\varepsilon}(L)=\pi t\varepsilon^2$. Next,  by
Theorem~\ref{th:3.3}(v), we have
\begin{equation}\label{e:4.11}
c_{HZ}^\circ({\rm D}_{\varepsilon}(E),
\bar\omega^t_{A})=c_{HZ}^\circ({\rm
D}_{1}(E),\omega^t_{A\varepsilon})
\end{equation}
for any $0<\varepsilon<\varepsilon_0(A, {\mathcal F})$ and
$0<t<t_0(A,\varepsilon)$.
 By the exact homotopy sequence of the fibration, we can easily
derive that the embedding $\varphi_{\mathcal P}:{\rm
D}_\delta(E)\to {\rm P}(E\oplus\C)$ induces an injective
homomorphism $\pi_1({\rm D}_\delta(E))\to\pi_1({\rm
P}(E\oplus\C))$ for any $0<\delta<1$. Applying (\ref{e:4.5}) to
$W_\delta:=\varphi_{\mathcal P}(Cl({\rm D}_\delta(E))))$, we get
\begin{equation}\label{e:4.11.1}
C_{HZ}^\circ(W_\delta,\Omega^t_{A\varepsilon})\le
C_{HZ}^\circ({\rm P}(E\oplus\C),\Omega^t_{A\varepsilon}; pt,
Z_0)\le\pi t\varepsilon^2
\end{equation}
for any $0<\varepsilon<\varepsilon_0({\mathcal F}, A)$ and $0<t<
t_0(A, \varepsilon)$. Here the second inequality comes from
(\ref{e:4.8}). Note that (\ref{e:4.1}) and
Theorem~\ref{th:3.3}(ii) imply
\begin{equation}\label{e:4.11.2}
c_{HZ}^\circ(Cl({\rm
D}_{\delta}(E)),\omega^t_{A\varepsilon})=C_{HZ}^\circ(Cl({\rm
D}_{\delta}(E)),\omega^t_{A\varepsilon})=C_{HZ}^\circ(W_\delta,\Omega^t_{A\varepsilon}),
\end{equation}
 for each $0<\varepsilon<\varepsilon_0({\mathcal F}, A)$ and
 $0<t<t_0(A,\varepsilon)$.
Setting $\delta\to 1$ and using the definition of $c_{HZ}^\circ$,
we obtain that
\begin{equation}\label{e:4.12}
c_{HZ}^\circ({\rm D}_{1}(E),\omega^t_{A\varepsilon})\le \pi
t\varepsilon^2
\end{equation}
for each $0<\varepsilon<\varepsilon_0({\mathcal F}, A)$ and
$0<t<t_0(A,\varepsilon)$. Now the desired (\ref{e:4.9}) follows
from (\ref{e:4.11}) and (\ref{e:4.12}) directly.

{\bf Step 2.} {\it Proving (\ref{e:4.10})}. Let $\dim B=2k$. Take
 a symplectic embedding \linebreak $\Upsilon:
(B^{2k}(r),\omega_{\rm std})\to (B,\beta)$, where $r={\mathcal
W}_G(B,\beta)$ is the Gromov symplectic radius. Since ${\rm
U}(E)|_{\Upsilon(B^{2k}(r))}$ can be trivialized, for any small
$\epsilon>0$,  we may choose a connection form $A'$ on ${\rm
U}(E)$ such that it is {\it flat} near
$\Delta:=Cl(\Upsilon(B^{2k}(r-\epsilon)))$. By
 Theorem~\ref{th:3.3}(vii), for any
 $0<\varepsilon\ll\min\{\varepsilon_0({\mathcal F}, A), \varepsilon_0({\mathcal F},
 A')\}\le 1$, $0<t<\min\{t_0(A,\varepsilon), t_0(A',\varepsilon)\}$, and $\delta\in (0,1]$,
we have a symplectomorphism
\begin{equation}\label{e:4.13}
\varphi^\delta_\varepsilon: ({\rm
 D}_\delta(E), \omega^t_{A\varepsilon})\to ({\rm
 D}_\delta(E), \omega^t_{A'\varepsilon}),
\end{equation}
which is also a bundle isomorphism. Note that $\Delta$ has the
same dimension as $B$. From the first paragraph of the proof of
Theorem~\ref{th:3.2}, it is easily observed that\linebreak
$\varepsilon_0^\circ(A'|_{Cl(\Delta)}, {\mathcal F}|_{{\rm
U}(E)|_{Cl(\Delta)}})$ can be taken as $\varepsilon_0(A',{\mathcal
F})$. Hence, Theorem~\ref{th:3.3}(vi) implies that $({\rm
 D}_1(E)|_\Delta, \omega^t_{A'|_\Delta\varepsilon})=({\rm
 D}_1(E)|_\Delta, \omega^t_{A'\varepsilon}|_\Delta)$ is a
 symplectic submanifold in $({\rm
 D}_1(E), \omega^t_{A'\varepsilon})$ for each
$0<\varepsilon<\varepsilon_0(A', {\mathcal F})$ and
$0<t<t_0(A',\varepsilon)$.   So for any symplectic
 capacity $c$ (including $c_{HZ}^\circ$), it holds that
\begin{equation}\label{e:4.14}
c({\rm
 D}_{1}(E)|_\Delta, \omega^t_{A'|_\Delta\varepsilon})\le
c({\rm
 D}_{1}(E), \omega^t_{A'\varepsilon}).
\end{equation}
  By Theorem~\ref{th:3.3}(v), the symplectomorphism
 \begin{equation}\label{e:4.14.1}
 \psi_\varepsilon: ({\rm
 D}_{1}(E), \omega^t_{A'\varepsilon})\to ({\rm
 D}_{\varepsilon}(E), \bar\omega^t_{A'}),\;(b,v)\mapsto (b,\varepsilon
 v)
 \end{equation}
maps $({\rm
 D}_{1}(E)|_\Delta, \omega^t_{A'|_\Delta\varepsilon})=({\rm
 D}_1(E)|_\Delta, \omega^t_{A'\varepsilon}|_\Delta)$ onto $({\rm
 D}_{\varepsilon}(E)|_\Delta, \bar\omega^t_{A'|_\Delta})$. We  get
\begin{equation}\label{e:4.15}
c({\rm
 D}_{1}(E)|_\Delta, \omega^t_{A'|_\Delta\varepsilon})=c({\rm
 D}_{\varepsilon}(E)|_\Delta, \bar\omega^t_{A'|_\Delta})
\end{equation}
for $0<\varepsilon<\varepsilon_0({\mathcal F}, A')$ and
$0<t<t_0(A',\varepsilon)$. Therefore (\ref{e:4.13}), (4.29), and
(\ref{e:4.15}) yields
\begin{eqnarray}\label{e:4.16}
c({\rm
 D}_{1}(E), \omega^t_{A\varepsilon})=c({\rm
 D}_{1}(E), \omega^t_{A'\varepsilon})\!\!\!\!\!\!\!\!\!\!&&\ge c({\rm
 D}_{1}(E)|_\Delta, \omega^t_{A'|_\Delta\varepsilon})\nonumber\\
&&= c({\rm
 D}_{\varepsilon}(E)|_\Delta, \bar\omega^t_{A'|_\Delta})
 \end{eqnarray}
for any
 $0<\varepsilon\ll\min\{\varepsilon_0({\mathcal F}, A), \varepsilon_0({\mathcal F},
 A')\}\le 1$ and $0<t<\min\{t_0(A,\varepsilon), t_0(A',\varepsilon)\}$.
Now Theorem~\ref{th:3.3}(viii) implies that for any
$0<\eta<\varepsilon_0(P, A_{\rm can})$ and $t>0$, the symplectic
manifold $({\rm
 D}_{\eta}(E)|_\Delta, \bar\omega^t_{A'|_\Delta})$ can be
 identified with
 \begin{equation}\label{e:4.16.1}
\bigl(\Delta\times B^{2n}(\eta), \beta\oplus(t\omega_{\rm
std}^{(n)})\bigr)\approx
 \bigr(Cl(B^{2k}(r-\epsilon))\times B^{2n}(\eta), \omega^{(k)}_{\rm
 std}\oplus(t\omega_{\rm std}^{(n)})\bigl).
 \end{equation}
 So for $0<t<{\mathcal W}_G(B,\beta)/\pi\eta^2$, it holds that
\begin{eqnarray}\label{e:4.17}
c({\rm
 D}_{\eta}(E)|_\Delta, \bar\omega^t_{A'|_\Delta})
\!\!\!\!\!\!\!\!\! &&=c\bigl(\Delta\times B^{2n}(\eta),
\beta\oplus(t\omega_{\rm
std}^{(n)})\bigr)\nonumber\\
&&=c\bigr(Cl(B^{2k}(r-\epsilon))\times B^{2n}(\eta),
\omega^{(k)}_{\rm
 std}\oplus(t\omega_{\rm std}^{(n)})\bigl)\nonumber\\
&&\ge{\mathcal W}_G\bigr(Cl(B^{2k}(r-\epsilon))\times
B^{2n}(\eta), \omega^{(k)}_{\rm
 std}\oplus(t\omega_{\rm std}^{(n)})\bigl)\nonumber\\
&&\ge \pi t\eta^2.
\end{eqnarray}
Taking $\eta=\varepsilon\ll\min\{\varepsilon_0({\mathcal F}, A),
\varepsilon_0({\mathcal F},
 A'), \varepsilon_0(P, A_{\rm can})\}$,  (\ref{e:4.16}) and
(\ref{e:4.17}) yield that
\begin{equation}\label{e:4.18}
c({\rm D}_{1}(E),\omega^t_{A\varepsilon})\ge \pi t\varepsilon^2
\end{equation}
for any $0<t<\min\{t_0(A',\varepsilon),
t_0(A,\varepsilon),{\mathcal W}_G(B,\beta)/\pi\varepsilon^2\}$.
Note that we can choose $\varepsilon>0$ so small that ${\mathcal
W}_G(B,\beta)/\pi\varepsilon^2>1$ and that we can also assume that
 $t_0(A',\varepsilon)$ and $t_0(A,\varepsilon)$ are no more than $1$.
The desired (\ref{e:4.10}) immediately follow from (\ref{e:4.9})
and (\ref{e:4.18}). \end{proof}

\noindent{\bf Proof of Theorem~\ref{th:1.3}.}\quad Clearly,
(\ref{e:1.6}) and (\ref{e:1.7}) directly  follow from
(\ref{e:4.9}) and (\ref{e:4.10}), respectively.
 To get (\ref{e:1.8}), for a
given $\epsilon>0$, we choose $t>0$ so small that $\pi t
\varepsilon^2<\epsilon$. Consider the (disk) bundle isomorphism
\begin{equation}\label{e:4.19}
\Theta_{t,\varepsilon}:{\rm
D}_{\sqrt{t}\varepsilon}(E)\to {\rm
D}_{\varepsilon}(E),\;(b,v)\mapsto (b,\frac{1}{\sqrt{t}} v).
\end{equation}
Then $\Theta_{t,\varepsilon}^\ast\bar\omega^t_{A}$ restricts to
$(\omega_{\rm std})_b$ on each fibre ${\rm
D}_{\sqrt{t}\varepsilon}(E)_b$, and the zero section $0_E$ is also
a symplectic submanifold in $({\rm
D}_{\sqrt{t}\varepsilon}(E),\Theta_{t,\varepsilon}^\ast\bar\omega^t_{A})$.
Hence, the zero section $0_E$ has the symplectic normal bundle
$(E,\bar\omega)$ in $({\rm
D}_{\sqrt{t}\varepsilon}(E),\Theta_{t,\varepsilon}^\ast\bar\omega^t_{A})$.
By the Weinstein's symplectic neighborhood theorem, there exists
$0<\eta\ll\varepsilon$ such that $({\rm
D}_{\sqrt{t}\eta}(E),\Theta_{t,\varepsilon}^\ast\bar\omega^t_{A})$
 and thus $({\rm
D}_{\eta}(E),\bar\omega^t_{A})\\
\subset ({\rm
D}_{\varepsilon}(E),\bar\omega^t_{A})$ is symplectomorphic to an
open neighborhood $W$ of $B$ in $M$. Then (\ref{e:1.8})  easily
follows from this and (\ref{e:1.7}).\hfill$\Box$\vspace{2mm}

\noindent{\bf Proof of Theorem~\ref{th:4.4}.}\quad  Since the
Gromov-Witten invariants are symplectic deformation invariants, it
suffices to prove that (\ref{e:4.7}) holds on $({\rm
P}(E\oplus\C),\Omega^t_{A\varepsilon})$  for sufficiently small
$t>0$. To this end, we need to construct a suitable
$\Omega^t_{A\varepsilon}$-tamed almost complex structure on ${\rm
P}(E\oplus\C)$. Using the standard complex structure $i$ and
symplectic structure on $\C$, we get a complex structure
$J_E\oplus i$ on $E\oplus\C$ and the corresponding Hermitian
structure. As in Remark~\ref{rm:4.3}, it in turn yields a complex
structure on each fibre ${\rm P}(E\oplus\C)_b$ of ${\rm
P}(E\oplus\C)$, still denoted by $J_{\rm can,b}$. One now has
obvious K\"ahler identification $\bigl({\rm P}(E\oplus\C)_b,
(\omega_{\rm FS})_b, J_{\rm can,b}\bigr)\equiv (\CP^n, \omega_{\rm
FS}, J_{\rm FS})$. From it and any compatible almost complex
structure $J_B\in {\mathcal J}(B,\beta)$ we can, as in
(\ref{e:4.6}), form an almost complex structure ${\bf J}$ on ${\rm
P}(E\oplus\C)$ which is $\Omega^t_{A\varepsilon}$-tamed for some
sufficiently small $0<t< t_0(A,\varepsilon)$. Moreover, the
projection ${\mathcal P}:({\rm P}(E\oplus\C), {\bf J})\to (B,
J_B)$ is an almost complex fibration with fibre $(\CP^n, J_{\rm
FS})$ by the final sentence above Remark~\ref{rm:4.3}. It is
well-known that $J_{\rm FS}\in{\mathcal J}_{\rm reg}(\CP^n,
\omega_{\rm FS})$, (cf. \cite[Prop.7.4.3]{McSa2}). Let
$[\CP^1]\subset H_2(\CP^n,\Z)$ denote the class of the line
$\CP^1\subset\CP^n$. It was proved in \cite{RT} that the
Gromov-Witten invariant $\Psi^{(\CP^n,\omega_{\rm FS})}_{[\CP^1],
0, 3}(pt; pt,pt, [\CP^{n-1}])=1$, where $pt$ always denotes the
class of a single point for different spaces, and $[\CP^{n-1}]\in
H_{2n-2}(\CP^n,\Z)$ is the class of the hyperplane
$\CP^{n-1}\subset\CP^n$.

Note that $Z_0$ (resp., $Z_\infty$) restricts to a single point
(resp., the hyperplane $\CP^{n-1}$) on each fiber ${\rm
P}(E\oplus\C)_b\equiv\CP^n$. If $L\in H_2({\rm P}(E\oplus\C),\Z)$
denotes the class of the line in the fiber, for the fibred almost
complex structure ${\bf J}$ just constructed  and a given point
$x\in {\rm P}(E\oplus\C)$, every ${\bf J}$-holomorphic sphere of
class $L$ in ${\rm P}(E\oplus\C)$ passing $x$ must sit in the
fiber ${\rm P}(E\oplus\C)_b$ at $b={\mathcal P}(x)\in B$.
Therefore, there exists such a unique curve passing $x$ and
generically intersecting with $Z_0$ and $Z_\infty$. As expected,
we arrive at (\ref{e:4.7}) for $({\rm
P}(E\oplus\C),\Omega^t_{A\varepsilon})$. Actually, by
\cite[Proposition 6.3.B]{Bi}, the ${\bf J}$ is regular for the
class $L$. So if  ${\mathcal M}({\rm P}(E\oplus\C), L, {\bf J})$
is the space of all ${\bf J}$-holomorphic spheres in ${\rm
P}(E\oplus\C)$ representing the class $L$, then ${\mathcal M}({\rm
P}(E\oplus\C), L, {\bf J})/{\rm PSL}(2,\R)$ is a compact smooth
manifold of dimension $\dim_{\R}B+ 2{\rm rank}_{\R}E-6$, and that
the dimension condition
\begin{eqnarray}\label{e:4.20}
\deg[pt]+\deg[Z_0]+\deg[Z_\infty]\!\!\!\!\!\!\!\!\!&&=\dim_{\R}{\rm
P}(E\oplus\C)+
2c_1(L)\nonumber\\
&&=\dim_{\R}B+ 2{\rm rank}_{\R}E+2
\end{eqnarray}
is satisfied. It easily follows that the Gromov-Witten invariant
(in the sense of \cite{RT}) \linebreak
$\Psi_{L,0,3}(pt;pt,[Z_0],[Z_\infty])=1$.

Since our \cite{Lu2} follows \cite{LiuT2}, the GW-invariant in
(\ref{e:4.7}) is one constructed by Liu-Tian in \cite{LiuT2} (also
see \cite{Lu6} for details). Let $\Psi^{\rm
vir}_{L,0,3}(pt;pt,[Z_0],[Z_\infty])$ denote this GW-invariant
constructed with the virtual moduli cycles. In the remainder of the
paper we shall show that it is also equal to $1$. (Namely, while
some GW-invariant can be simultaneously defined in the methods in
[14, 34], they agree.) Instead of using the method in \cite{Lu5}, we
use the method of proof of \cite[Proposition 7.6]{Lu2} to prove that
the GW-invariant $\Psi^{\rm vir}_{L,0,3}(pt;pt,[Z_0],[Z_\infty])$ in
the sense of the general definition is also equal to $1$. Let
$\overline{\mathcal M}_{0,3}({\rm P}(E\oplus\C), L, {\bf J})$ be the
space of all equivalence classes of all $3$-pointed stable ${\bf
J}$-maps of genus zero and of class $L$ in ${\rm P}(E\oplus\C)$.
Then each stable map $[{\bf f}]\in\overline{\mathcal M}_{0,3}({\rm
P}(E\oplus\C), L, {\bf J})$ has an image set contained in a single
fiber of ${\rm P}(E\oplus\C)$ since the image set is connected.
Assume that ${\rm Im}({\bf f})\subset {\rm P}(E\oplus\C)_b$. This
$[{\bf f}]$ may naturally be viewed as an element of
$\overline{\mathcal M}_{0,3}({\rm P}(E\oplus\C)_b, L_b, J_{\rm
can,b})$, where $L_b$ is the class of the line in ${\rm
P}(E\oplus\C)_b\equiv\CP^n$. Since $\Omega^t_{A\varepsilon}|_{{\rm
P}(E\oplus\C)_b}=t\varepsilon^2(\omega_{\rm FS})_b$ and $L_b$ is
indecomposable with respect to $(\omega_{\rm FS})_b$, by the proof
of \cite[Proposition 7.6]{Lu2}, it is easily seen that
$\overline{\mathcal M}_{0,3}({\rm P}(E\oplus\C)_b, L_b,  J_{\rm
can,b})$ is a stratified smooth compact manifold (using the
regularity of $J_{\rm can,b}$). More precisely,
\begin{equation}\label{e:4.21}
\overline{\mathcal M}_{0,3}({\rm P}(E\oplus\C)_b, L_b,  J_{\rm
can,b})=\bigcup^4_{i=1}\overline{\mathcal M}_{0,3}({\rm
P}(E\oplus\C)_b, L_b,  J_{\rm can,b})_i,
\end{equation}
where each stratum $\overline{\mathcal M}_{0,3}({\rm
P}(E\oplus\C)_b, L_b, J_{\rm can,b})_i$ is a smooth manifold, and
\begin{equation}\label{e:4.22}
\left.\begin{array}{ll}
 &\dim \overline{\mathcal M}_{0,3}({\rm
P}(E\oplus\C)_b, L_b, J_{\rm
can,b})_1=2n + 2(n+1),\\
&\dim \overline{\mathcal M}_{0,3}({\rm P}(E\oplus\C)_b, L_b,
J_{\rm
can,b})_2=2n + 2(n+1)-4,\\
&\dim \overline{\mathcal M}_{0,3}({\rm P}(E\oplus\C)_b, L_b,
J_{\rm
can,b})_3=2n + 2(n+1)-6,\\
&\dim \overline{\mathcal M}_{0,3}({\rm P}(E\oplus\C)_b, L_b,
J_{\rm can,b})_4=2n + 2(n+1)-6.\end{array}\right.
\end{equation}
So $\overline{\mathcal M}_{0,3}({\rm P}(E\oplus\C), L, {\bf
J})=\bigcup_{b\in B}\overline{\mathcal M}_{0,3}({\rm
P}(E\oplus\C)_b, L_b, J_{\rm can,b})$ is also a stratified smooth
compact manifold (with correct dimension because of the regularity
of ${\bf J}$). That is,
\begin{equation}\label{e:4.23}
\overline{\mathcal M}_{0,3}({\rm P}(E\oplus\C), L, {\bf
J})=\bigcup^4_{i=1}\overline{\mathcal M}_{0,3}({\rm P}(E\oplus\C),
L, {\bf J})_i,
\end{equation}
where $\overline{\mathcal M}_{0,3}({\rm P}(E\oplus\C), L, {\bf
J})_i=\bigcup_{b\in B}\overline{\mathcal M}_{0,3}({\rm
P}(E\oplus\C)_b, L_b, J_{\rm can,b})_i$, $i=1,2,3,4$. Since each
stable map  $[{\bf f}]\in\overline{\mathcal M}_{0,3}({\rm
P}(E\oplus\C), L, {\bf J})$ has
 no free components, we  have the following claim.

 \begin{claim}\label{cl:4.6}
 A virtual moduli cycle of $\overline{\mathcal
M}_{0,3}({\rm P}(E\oplus\C), L, {\bf J})$ can be taken as
\begin{equation}\label{e:4.24}
\overline{\mathcal M}_{0,3}({\rm P}(E\oplus\C), L, {\bf J})\to
{\mathcal B}_{0,3,L}^{{\rm P}(E\oplus\C)},\;[{\bf f}]\mapsto [{\bf
f}].
\end{equation}
\end{claim}

 Once it is proved, then almost repeating the proof of \cite[Proposition
 7.6]{Lu2}, we can obtain the desired (\ref{e:4.7}):
\begin{equation}\label{e:4.25}
\Psi^{\rm
vir}_{L,0,3}(pt;pt,[Z_0],[Z_\infty])=\Psi_{L,0,3}(pt;pt,[Z_0],[Z_\infty])=1.
\end{equation}

In order to prove Claim~\ref{cl:4.6}, as in the proof of
\cite[Proposition 7.6]{Lu2}, each \linebreak $[{\bf
f}]\in\overline{\mathcal M}_{0,3}({\rm P}(E\oplus\C), L, {\bf J})$
must be one of the following four cases.
\begin{itemize}
\item[(i)] The domain $\Sigma=\CP^1$, $z_i$, $i=1,2,3$, are three
distinct marked points on $\Sigma$, and $f:\Sigma\to {\rm
P}(E\oplus\C)$ is a $J$-holomorphic map of class $L$.

\item[(ii)] The domain $\Sigma$ has exactly two components
$\Sigma_1=\CP^1$ and $\Sigma_2=\CP^1$ which have a unique
intersecting point. $f|_{\Sigma_1}$ is nonconstant and $\Sigma_1$
only contains a marked point. $f|_{\Sigma_2}$ is constant and
$\Sigma_2$ contains two marked points.

\item[(iii)] The domain $\Sigma$ has exactly two components
$\Sigma_1=\CP^1$ and $\Sigma_2=\CP^1$ which have a unique
intersecting point. $f|_{\Sigma_1}$ is nonconstant and $\Sigma_1$
contains no  marked point. $f|_{\Sigma_2}$ is constant and
$\Sigma_2$ contains three marked points.

\item[(iv)] The domain $\Sigma$ has exactly three components
$\Sigma_1=\CP^1$, $\Sigma_2=\CP^1$, and $\Sigma_3=\CP^1$.
$\Sigma_1$ and $\Sigma_2$ (resp., $\Sigma_2$ and $\Sigma_3$) has
only an intersecting point, and $\Sigma_1$ and $\Sigma_3$ has no
intersecting point.  $f|_{\Sigma_1}$ is nonconstant and $\Sigma_1$
contains no a marked point. $f|_{\Sigma_2}$ is constant and
$\Sigma_2$ contains a marked point. $f|_{\Sigma_3}$ is constant
and $\Sigma_3$ contains two marked points.
\end{itemize}
Note that in each case, the nonconstant $f|_{\Sigma_1}$ is simple,
and thus somewhere injective. It follows that the automorphism
group $Aut({\bf f})$ of ${\bf f}$ is trivial though $[{\bf f}]$
might contain  many representatives. Let ${\mathcal
B}_{L,0,3}^{{\rm P}(E\oplus\C)}$ be the set of equivalence classes
of all $3$-pointed stable $L^{k,p}$-maps in ${\rm P}(E\oplus\C)$
which represent class $L$, have genus $0$ and domains that belong
to the four types above. For $[{\bf f}]\in\overline{\mathcal
M}_{0,3}({\rm P}(E\oplus\C), L, {\bf J})$, let $\widetilde {\bf
U}_\delta({\bf f}, {\rm H})$ be the local uniformizer near $[{\bf
f}]\in {\mathcal B}_{L,0,3}^{{\rm P}(E\oplus\C)}$ as constructed
in \cite[Section 2]{LiuT1} (also see \cite[Section 2]{Lu6} for
details). Since $Aut({\bf f})$ is trivial, \cite[Lemma 2.6]{LiuT1}
showed that the natural projection $\widetilde {\bf U}_\delta({\bf
f}, {\rm H})\to {\mathcal B}_{L,0,3}^{{\rm P}(E\oplus\C)}$ gives
rise to a homeomorphism to an open neighborhood of $[{\bf f}]$ in
${\mathcal B}_{L,0,3}^{{\rm P}(E\oplus\C)}$ (see \cite[Section
4.2.1]{Mc} for more explanations). This suggests that some open
neighborhood ${\mathcal W}$ of $\overline{\mathcal M}_{0,3}({\rm
P}(E\oplus\C), L, {\bf J})$ in ${\mathcal B}_{L,0,3}^{{\rm
P}(E\oplus\C)}$ might carry a stratified Banach manifold
structure. Regardless of these, we still adopt the Liu-Tian
construction method in \cite{LiuT1, LiuT2} to get a system of
bundles
\begin{equation}\label{e:4.26}
(\widetilde{\mathcal E}^\Gamma, \widetilde W^\Gamma)=
\bigl\{(\widetilde E_I^{\Gamma_I}, \widetilde W_I^{\Gamma_I}),
\tilde\pi_I, \widetilde\Pi_I, \Gamma_I, \tilde\pi^I_J,
\widetilde\Pi^I_J,\,\lambda^I_J\,\bigm|\, J\subset I\in{\mathcal
N}\bigr\}
\end{equation}
(cf. \cite{Lu6}). However, so far all $\widetilde
E_I^{\Gamma_I}\to\widetilde W_I^{\Gamma_I}$ are stratified Banach
bundles on stratified Banach manifolds, all groups $\Gamma_I$ and
homomorphisms $\lambda^I_J$ are trivial, and all projections
$\tilde\pi_I$, $\widetilde\Pi_I$,  $\tilde\pi^I_J$ and
$\widetilde\Pi^I_J$ becomes open embeddings. Therefore, it is not
needed to make furthermore desingularization for the system of
bundles in (4.43) as done in \cite{LiuT3} (also see \cite{Lu6} for
a detailed description). We only need to renormalize
$(\widetilde{\mathcal E}^\Gamma, \widetilde W^\Gamma)$ to get a
new system of stratified smooth Banach bundles
\begin{equation}\label{e:4.27} (\widetilde
{\mathcal E}, \widetilde V)=\bigl \{(\widetilde E_I, \widetilde
V_I), \tilde\pi_I,
 \tilde\pi^I_J, \widetilde\Pi_I, \widetilde\Pi^I_J, \tilde p_I,\Gamma_I\,\bigm|\, J\subset I\in{\mathcal
 N}\bigr\}
\end{equation}
 as in \cite{LiuT1, LiuT2}. Then by restrictions, we get a
system of stratified smooth Banach bundles as in
\cite[(3.74)]{Lu6}:
\begin{equation}\label{e:4.28} (\widetilde {\mathcal
E}^\ast, \widetilde V^\ast)=\bigl \{(\widetilde E_I^\ast,
\widetilde V_I^\ast), \tilde\pi_I,
 \tilde\pi^I_J, \widetilde\Pi_I, \widetilde\Pi^I_J, \tilde p_I,\Gamma_I\,\bigm|\,
 J\subset I\in{\mathcal N}\bigr\}.
\end{equation}
Now following the idea of the proof of \cite[Theorem 4.1]{LiuT1}
(see \cite[(3.75) and (3.76)]{Lu6}), we have the obvious pullback
stratified smooth Banach bundle system
\begin{equation}\label{e:4.29} \bigl({\bf P}_1^\ast\widetilde
{\mathcal E}^\ast, \widetilde V^\ast\times {\bf
B}_\eta({\R}^q)\bigr) =\bigl \{({\bf P}_1^\ast\widetilde E_I^\ast,
\widetilde V_I^\ast\times {\bf B}_\eta({\R}^q)), \pi_I, \pi^I_J,
\Pi_I,\Pi^I_J, p_I,\Gamma_I\,\bigm|\, J\subset I\in{\mathcal
N}\bigr\},
\end{equation}
and its  global section $\Psi=\{\Psi_I\,|\, I\in{\mathcal N}\}$,
\begin{equation}\label{e:4.30}
\Psi_I: \widetilde V_I^\ast\times {\bf B}_\eta({\R}^q)\to {\bf
P}_1^\ast\widetilde E_I^\ast,\; (\tilde x_I, {\bf t})\mapsto
(\tilde\partial_{\bf J})_I(\tilde x_I)+
 \sum^{n_3}_{i=1}\sum^{q_i}_{j=1}t_{ij}(\tilde s_{ij})_I(\tilde
 x_I),
\end{equation}
where ${\bf t}=\{t_{ij}| 1\le j\le q_i,\,1\le i\le n_3\}\in
{\R}^q$.
  Clearly, $\Psi_I(\tilde x_I, 0)=0$ for any zero $\tilde x_I$ of
 $(\tilde\partial_{\bf J})_I$ in $\widetilde V_I$. By \cite[Theorem 3.14 and
 Corollary 3.15]{Lu6}, we get a small $\eta>0$ and a residual
subset ${\bf B}^{res}_\eta({\R}^q)\subset {\bf B}_\eta({\R}^q)$
such that for each ${\bf t}\in {\bf B}^{res}_\eta({\R}^q)$, the
global section
 $\Psi^{({\bf t})}=\{\Psi^{({\bf t})}_I\,|\, I\in{\mathcal
N}\}$ of the bundle system $\bigl(\widetilde {\mathcal E}^\ast,
\widetilde V^\ast\bigr)$ is transversal to the zero section, where
$ \Psi^{({\bf t})}_I: \widetilde V_I^\ast\to \widetilde
E_I^\ast,\; \tilde x_I \mapsto \Psi_I(\tilde x_I,{\bf t})$. So the
set $\widetilde{\mathcal M}^{\bf t}_I:= (\Psi^{({\bf
t})}_I)^{-1}(0)$ is a stratified smooth Banach manifold of
dimension $\dim B+ 4n+ 2$. It also holds that
\begin{itemize}
\item[(A)] the stratified Banach manifold $\widetilde{\mathcal
M}_I^{{\bf t}}$ has no strata of codimension odd, and each stratum
of $\widetilde{\mathcal M}^{\bf t}_I$ of codimension $r$ is
exactly the intersection of $\widetilde{\mathcal M}^{\bf t}_I$ and
the stratum of $\widetilde V_I^\ast$ of codimension $r$ for
$r=0,\cdots, \dim B+ 4n+ 2$;

\item[(B)] the family $\widetilde{\mathcal M}^{\bf
t}=\bigl\{\widetilde{\mathcal M}^{\bf t}_I\,|\, I\in{\mathcal
N}\bigr\}$ is compatible in the sense that for any $J\subset
I\in{\mathcal N}$,
\begin{equation}\label{e:4.31}
 \tilde\pi^I_J: (\tilde \pi^I_J)^{-1}\bigl(\widetilde{\mathcal M}^{\bf
t}_I\bigr)\to {\rm Im}(\tilde\pi^I_J)\subset\widetilde{\mathcal
M}^{\bf t}_J
 \end{equation}
 is a continuous and stratified smooth open embedding;

\item[(C)]  for each $I\in{\mathcal N}$ and any two ${\bf t},\,
{\bf t}'\in {\bf B}^{res}_\eta({\R}^q)$, the cornered stratified
Banach manifolds $\widetilde{\mathcal M}_I^{{\bf t}}$ and
$\widetilde{\mathcal M}_I^{{\bf t}'}$ are cobordant, and thus maps
$\tilde\pi_I: \widetilde{\mathcal M}_I^{{\bf t}}\to {\mathcal W}$
and $\tilde\pi_I: \widetilde{\mathcal M}_I^{{\bf t}'}\to {\mathcal
W}$ are also cobordant.
\end{itemize}
Since $|\Gamma_I|=1$ for any
$I\in{\mathcal N}$, the formal summations
\begin{equation}\label{e:4.32}
{\mathcal C}^{\bf t}:=\sum_{I\in{\mathcal
 N}}\{\tilde\pi_I:\widetilde{\mathcal M}^{\bf t}_I\to{\mathcal
 W}\}\quad\forall{\bf t}\in {\bf B}_\varepsilon^{res}({\R}^q),
 \end{equation}
a family of  cobordant singular cycles in
 ${\mathcal W}$,  are virtual moduli cycles in ${\mathcal
 W}$ constructed by Liu-Tian method (cf. [12, 13, 14, 21]).
 As explained  \cite[page 65]{LiuT1}, the summation
precisely means that on the overlap of two pieces of ${\mathcal
C}^{\bf t}$, we only count them once. Since all $\tilde\pi_I$ and
$\tilde\pi^I_J$ are stratified smooth open embeddings,
$\{(\tilde\pi_I,{\mathcal M}^{\bf t}_I)\,|\, I\in{\mathcal
 N}\}$ is actually a compatible coordinate chart cover of
 a compact, stratified smooth Banach manifold
 $\cup_{I\in{\mathcal N}}\tilde\pi_I({\mathcal M}^{\bf
 t}_I)\subset{\mathcal W}$ of dimension $\dim B+ 4n+2$ and without
 strata of codimension one.

Note that each $(\tilde\partial_{\bf J})_I$ is essentially the
Cauchy-Riemann operator $\bar\partial_{\bf J}$ and that ${\bf J}$
is regular with respect to the class $L$. It is not hard to prove
that ${\bf 0}\in {\bf B}_\varepsilon^{res}({\R}^q)$ and
$\cup_{I\in{\mathcal
 N}}\tilde\pi_I(\widetilde{\mathcal M}^{\bf 0}_I)=\overline{\mathcal M}_{0,3}({\rm P}(E\oplus\C), L, {\bf
 J})$. This proved Claim~\ref{cl:4.6}.
\hfill$\Box$

\begin{remark}\label{rm:4.7}{\rm As suggested in the above proof,
so far  ${\mathcal E}\to {\mathcal W}$ might be a stratified
Banach bundle on a stratified Banach manifold.
 The Cauchy-Riemann operator
$\bar\partial_{\bf J}$ is a Fredholm section (restricting each
stratum) and has $\overline{\mathcal M}_{0,3}({\rm P}(E\oplus\C),
L, {\bf J})$ as zero set of it. The original arguments of
transversality and gluing may yield  finitely many  continuous
stratified smooth sections $s_i:{\mathcal W}\to {\mathcal E}$,
$i=1,\cdots, m$ such that the section
\begin{equation}\label{e:4.33} \Phi:{\mathcal W}\times {\rm
B}_\eta(\R^m)\to\Pi_1^\ast{\mathcal E},\;(\tau, {\bf
t})\mapsto\bar\partial_{\bf J}\tau+ t_1s_1(\tau)+\cdots +
t_ms_m(\tau)
\end{equation}
is transversal to the zero section for $\eta>0$ small enough.
Consequently, for generic small ${\bf t}\in {\rm B}_\eta(\R^m)$,
the section $\Phi_{\bf t}:{\mathcal W}\to {\mathcal
E},\;\tau\mapsto \bar\partial_{\bf J}\tau+ t_1s_1(\tau)+\cdots +
t_ms_m(\tau)$ is transversal to the zero section. In particular,
since ${\bf J}$ is
 regular, the section $\Phi_0=\bar\partial_{\bf J}$ is transversal to
the zero section. So we do not need to use Liu-Tian method as
above and can construct a desired perturbation cycle of
$\overline{\mathcal M}_{0,3}({\rm P}(E\oplus\C), L, {\bf J})$,
 which is cobordant to virtual moduli
cycles constructed by Liu-Tian method.  These are, in detail,
explained and developed in general abstract settings in
\cite{LuT}. }
\end{remark}

\vspace{2mm}

\noindent{\bf Acknowledgements}.\quad This work was partially
supported by the NNSF 10371007 of China and the Program for New
Century Excellent Talents of the Education Ministry of China. This
work began during the author's visit at the Institut des Hautes
\'Etudes Scientifiques (IHES) in Winter 2004 and was completed
during the author's visit to the Abdus Salam International Centre
for Theoretical Physics in
 Summer 2005. The author cordially thanks  Professors J. P.
Bourguignon and D.T. L\^e  for their invitations and for the
financial support and hospitality of their institutions, and
Professor C. Viterbo for his invitation to join Symplectic Geometry
Seminar in \'Ecole Polytechnique during the author's visit to IHES
in Winter 2004.  The author also thanks Paul Biran, Leonardo
Macarini, Felix Schlenk and Ramadas Ramakrishnan for some related
discussions. The author would like to thank Professor Leonid
Polterovich for mailing  me his lovely book \cite{Po3} and for very
valuable improvement suggestions, and he would like to thank Ely
Kerman for sending me his preprint. Thanks are also due to the
anonymous referee, whose comments and corrections improved the
exposition.

\end{document}